\documentclass[10pt]{article}
\usepackage{amsfonts,color}
\usepackage{amssymb}
\usepackage{footnote}
\usepackage{mathtools}
\usepackage{amsthm,wasysym}
\usepackage[english]{babel}
\usepackage{bbm}
\usepackage{colonequals}
\usepackage{setspace}
\usepackage{titlesec}

\usepackage{graphicx}
\usepackage{amsfonts}
\usepackage{hyperref}
\usepackage{subcaption}
\usepackage{amsmath}
\usepackage{nicefrac}
\usepackage{gensymb}
\usepackage{enumerate}
\usepackage[nameinlink,capitalize]{cleveref}
\usepackage{chngcntr}
\usepackage{float}

\usepackage{cite}
\usepackage[nottoc,numbib]{tocbibind}

\usepackage{pgf,tikz,pgfplots}
\pgfplotsset{compat=1.14}
\usepackage{mathrsfs}
\usetikzlibrary{arrows}
\usetikzlibrary{arrows.meta}

\allowdisplaybreaks[1]

\makeindex

\newtheoremstyle{break}
{\topsep}{\topsep}%
{\itshape}{}%
{\bfseries}{}%
{\newline}{}%
\theoremstyle{break}

\newtheorem{theorem}{Theorem}

\newtheorem{remark}{Remark}

\newtheorem{definition}{Definition}

\counterwithin{theorem}{section}
\counterwithin{lemma}{section}
\counterwithin{remark}{section}
\counterwithin{observation}{section}
\counterwithin{definition}{section}

\counterwithin{figure}{section}
\counterwithin{equation}{section}

\makeatletter
\let\@fnsymbol\@arabic
\makeatother

\crefname{Problem}{Problem.}{Problem.}

\DeclareMathOperator{\at}{\bigg\vert}
\newcommand{\vb}[1]{\mathbf{#1}}
\newcommand{\bm}[1]{\boldsymbol{#1}}

\DeclareMathOperator{\sym}{\mathrm{sym}}
\DeclareMathOperator{\spa}{\mathrm{span}}
\DeclareMathOperator{\sph}{\mathrm{sph}}
\DeclareMathOperator{\dev}{\mathrm{dev}}
\DeclareMathOperator{\skw}{\mathrm{skew}}
\DeclareMathOperator{\Anti}{\mathrm{Anti}}
\DeclareMathOperator{\tr}{\mathrm{tr}}

\newcommand{\jump}[1]{\ensuremath{[\![#1]\!]} }
\newcommand{\one}{\bm{\mathbbm{1}}}
\newcommand{\con}[2]{\langle {#1} , \, {#2} \rangle}
\newcommand{\norm}[1]{\| {#1} \|}

\newcommand{\dd}{\mathrm{d}}
\newcommand{\D}{\mathrm{D}}
\DeclareMathOperator{\di}{\mathrm{div}}
\DeclareMathOperator{\Di}{\mathrm{Div}}

\DeclareMathOperator{\curl}{\mathrm{curl}}
\DeclareMathOperator{\Curl}{\mathrm{Curl}}

\newcommand{\so}{\mathfrak{so}}
\renewcommand{\sl}{\mathfrak{sl}}

\newcommand{\Sph}{\R \cdot \one}

\newcommand{\Nedtwo}{\mathcal{N}_{II}}
\newcommand{\CG}{\mathcal{CG}}
\newcommand{\DG}{\mathcal{DG}}

\newcommand{\Po}{\mathit{P}}

\newcommand{\Le}{{\mathit{L}^2}}
\newcommand{\Lez}{{\mathit{L}^2_0}}

\newcommand{\Hone}{\mathit{H}^1}
\newcommand{\Honez}{\mathit{H}_0^1}

\newcommand{\Hc}[1]{\mathit{H}(\mathrm{curl}{#1})}
\newcommand{\HC}[1]{\mathit{H}(\mathrm{Curl}{#1})}

\newcommand{\HdD}[1]{\mathit{H}(\di\mathrm{Div}{#1})}
\newcommand{\HsC}[1]{\mathit{H}(\mathrm{sym}\,\mathrm{Curl}{#1})}

\newcommand{\HsCd}[1]{\mathit{H}^{\mathrm{dev}}(\mathrm{sym}\,\mathrm{Curl}{#1})}
\newcommand{\HCd}[1]{\mathit{H}^{\mathrm{dev}}(\mathrm{Curl}{#1})}
\newcommand{\HCdz}[1]{\mathit{H}_0^{\mathrm{dev}}(\mathrm{Curl}{#1})}
\newcommand{\HCv}[1]{\mathit{H}^{\mathrm{sph}}(\mathrm{Curl}{#1})}
\newcommand{\HsCv}[1]{\mathit{H}^{\mathrm{sph}}(\mathrm{sym}\,\mathrm{Curl}{#1})}

\newcommand{\elem}{T}
\newcommand{\body}{V}
\newcommand{\surf}{A}
\newcommand{\curv}{s}

\newcommand{\T}{\mathit{T}}

\newcommand{\R}{\mathbb{R}}

\newcommand{\C}{\mathit{C}}
\newcommand{\X}{\mathit{X}}
\renewcommand{\H}{\mathit{H}}

\newcommand{\Y}{\mathcal{Y}}
\newcommand{\Z}{\mathit{Z}}
\newcommand{\W}{\mathit{W}}

\newcommand{\lame}{\lambda_{\mathrm{e}}}
\newcommand{\lammi}{\lambda_{\mathrm{m}}}
\newcommand{\lamma}{\lambda_{\mathrm{M}}}
\newcommand{\mue}{\mu_{\mathrm{e}}}
\newcommand{\muc}{\mu_{\mathrm{c}}}
\newcommand{\mumi}{\mu_{\mathrm{m}}}
\newcommand{\muma}{\mu_{\mathrm{M}}}
\newcommand{\Lc}{L_\mathrm{c}}
\newcommand{\Ce}{\mathbb{C}_{\mathrm{e}}}
\newcommand{\Cc}{\mathbb{C}_{\mathrm{c}}}
\newcommand{\Cmic}{\mathbb{C}_{\mathrm{m}}}
\newcommand{\Cmac}{\mathbb{C}_{\mathrm{M}}}

\newcommand{\ud}{\vb{u}}
\newcommand{\Pm}{\bm{P}}

\newcommand{\sphP}{\bm{\mathfrak{I}}}
\newcommand{\devP}{\bm{\mathfrak{D}}}
\newcommand{\hP}{\bm{\mathfrak{h}}}
\newcommand{\pdevP}{\bm{\mathfrak{d}}}

\title{Novel $H^\mathrm{dev}(\mathrm{Curl})$-conforming elements on regular triangulations and Clough--Tocher splits for the planar relaxed micromorphic model}

\author{\normalsize{Adam Sky}\thanks{Corresponding author: Adam Sky, Institute of Computational Engineering and Sciences, Department of Engineering, Faculty of Science, Technology and Medicine, University of Luxembourg, 6, Avenue de la Fonte, L-4362 Esch-sur-Alzette, Luxembourg, email: adam.sky@uni.lu}
	, \quad
	\normalsize{Michael Neunteufel}\thanks{Michael Neunteufel, Department of Mathematics and Statistics, Portland State University, Portland, OR, United States, email: mneunteu@pdx.edu} 
	, \quad
    \normalsize{Peter Lewintan}\thanks{Peter Lewintan, Karlsruhe Institute of Technology, Englerstrasse 2, 76131 Karlsruhe, Germany, email: peter.lewintan@kit.edu}
    , \quad 
    \normalsize{Panos Gourgiotis}\thanks{Panos Gourgiotis, Mechanics Division, SAMPS, National Technical University of Athens, Zographou, GR-15773, Greece, email: pgourgiotis@mail.ntua.gr}
    , 
    \\
    \normalsize{Andreas Zilian}\thanks{Andreas Zilian, Institute of Computational Engineering and Sciences, Department of Engineering, Faculty of Science, Technology and Medicine, University of Luxembourg, 6, Avenue de la Fonte, L-4362 Esch-sur-Alzette, Luxembourg, email: andreas.zilian@uni.lu}
	\quad
	\normalsize{and} \quad
	\normalsize{Patrizio Neff}\thanks{Patrizio Neff, Chair for Nonlinear 
		Analysis and Modelling, Faculty of Mathematics, Universit\"{a}t Duisburg-Essen,
		Thea-Leymann Str. 9, 45127 Essen, Germany, email: patrizio.neff@uni-due.de}
}

\begin{document}

\maketitle

\begin{abstract}
In this work we present a consistent reduction of the relaxed micromorphic model to its corresponding two-dimensional planar model, such that its capacity to capture discontinuous dilatation fields is preserved. As a direct consequence of our approach, new conforming finite elements for $\HCd{,\surf}$ become necessary. We present two novel $\HCd{,\surf}$-conforming finite element spaces, of which one is a macro element based on Clough--Tocher splits, as well as primal and mixed variational formulations of the planar relaxed micromorphic model. Finally, we demonstrate the effectiveness of our approach with two numerical examples. 
\\
\vspace*{0.25cm}

\end{abstract}

\section{Introduction}
The accurate modelling of materials with complex micro-structures is a significant challenge in computational material science \cite{SARHIL2024116944,Sarhil2023}. Capturing the detailed kinematics of these materials requires precise geometric modelling of the entire domain, encompassing unit-cell geometries in metamaterials \cite{Voss,RIZZI2024105269,ALBERDI2021104540} or the various hole-shapes in porous media \cite{Seibert,Seibert2024,Boon2023}. The level of detail directly affects the resolution of the discretisation in finite element simulations, which can lead to increased computational effort \cite{Danesh}. Consequently, various methods are contemporarily applied to actively circumvent or reduce these computational costs, such as multi-scale simulations \cite{LIU2021114161,LIU2022105018}, or enriched continuum models \cite{Neff2014,Trinh2012}.    

The relaxed micromorphic model \cite{Neff2014,SKY2022115298,sky_hybrid_2021,GOURGIOTIS2024112700,Schroder2022,KNEES2023126806} represents one approach in the field of continuum mechanics that incorporates the micro-structural behaviour of materials directly into its kinematical description. This is achieved by taking possible distortions of the solution into account directly in the mathematical model via the introduction of a so called microdistortion field $\Pm$.
Thus, the model captures averaged micro-kinematics without having to explicitly resolve their geometries.
The model distinguishes itself from other micromorphic models \cite{Forest,Forest2006,ALAVI2021104278,Alavi2023,Li2020,boon2024mixed,Ghiba} by taking only the Curl of the microdistortion field $\Curl \Pm$ into consideration in its energy functional as opposed to the full gradient $\D \Pm$ in the classical micromorphic model by Eringen and Mindlin \cite{Eringen1999,Mindlin1964}.
This leads to many advantages: the micro-dislocation field $\Curl\Pm$ remains a second order tensor, the amount of new material coefficients is reduced, a physically consistent boundary condition for the microdistortion field $\Pm$ is made possible \cite{d’Agostino2022}, and the corresponding space of the microdistortion field is enlarged
\begin{align}
    &\Pm \in  \underbrace{\{ \Pm \in [\Le(\body)]^{3 \times 3} \; | \; \Curl\Pm \in [\Le(\body)]^{3 \times 3} \}}_{\HC{,\body}} \supset  \underbrace{\{\Pm \in [\Le(\body)]^{3 \times 3} \; | \; \D\Pm \in [\Le(\body)]^{3 \times 3 \times 3}\}}_{[\Hone(\body)]^{3 \times 3}} \, , && \body \subset \R^3 \, ,
\end{align}
where $\D \bm{P}: \body \subset \R^3 \to \R^{3 \times 3 \times 3}$ is a third-order tensor given by $\bm{P}\otimes \nabla = \bm{P}_{,i}\otimes \vb{e}_i$, and $\Le(\body)$ denotes the classical Lebesgue space of square integrable functions over $\body\subset \R^3$
\begin{align}
    \Le(\body) &=  \{ u: \body \to \mathbb{R} \, | \,  \|u\|_{\Le}^2 < \infty \} \, , & \|u \|_{\Le}^2 &= \int_{\body} \|u\|^2 \, \dd \body \, .
\end{align}
A possibility to further relax the assumptions on the regularity of the microdistortion by taking $\sym\Curl \Pm$ instead of $\Curl \Pm$ was envisaged in \cite{Lewintan2021}. In fact, this relaxation leads to the even larger space 
\begin{align}
    &\Pm \in  \underbrace{\{ \Pm \in [\Le(\body)]^{3 \times 3} \; | \; \sym\Curl\Pm \in [\Le(\body)]^{3 \times 3} \}}_{\HsC{,\body}} \supset \HC{,\body} \, ,
\end{align}
for the microdistortion field $\Pm$, while simultaneously maintaining well-posedness due to the generalised incompatible Korn-type inequalities \cite{Lewintan2021,LewintanInc,LewintanInc2,NEFF20151267,Neff2012,Gmeineder1,Gmeineder2}. As a result, $\HsC{,\body}$-conforming finite elements were developed in \cite{SKYNOVEL} specifically for the relaxed micromorphic model. 

In \cite{SKYNOVEL,SkyOn,Lewintan2021} we observed that $\Pm \in \HsC{,\body}$ allows for jumping dilatation fields. However, this property is not inherited in the dimensional reduction of the model to plane strain. This is because the operator $\sym\Curl$ is meaningless in two dimensions, since the two-dimensional $\Curl$-operator maps a second order tensor to a vector. In this work we reconsider this feature by re-writing the three-dimensional model in a slightly different manner, such that the dimensional reduction maintains the discontinuous dilatation property. As a result, we require novel finite elements for computations of the two-dimensional model. Namely, we require conforming subspaces for
\begin{align}
    &\HCd{,\surf} = \{ \devP \in \HC{,\surf}  \; | \; \tr \devP = 0 \} \, , 
\end{align}
where 
\begin{align}
    \HC{,\surf} = \{ \Pm \in [\Le(\surf)]^{2 \times 2} \; | \; \Curl \Pm \in [\Le(\surf)]^2 \} \, , && \surf \subset \R^2 \, .
\end{align}
We address their construction with two possibilities, of which the second one relies on a macro element approach \cite{gopalakrishnan2024johnsonmercier,Christiansen2024,Scott} on Clough--Tocher splits \cite{Philippe}.

This work is organised as follows: We shortly introduce the full three-dimensional relaxed micromorphic model. The model is subsequently re-written and reduced for two-dimensional domains, where we also explore its behaviour for the limits of the characteristic length-scale parameter $\Lc$. We then prove well-posedness of the reduced model for both a primal and mixed forms. After dealing with the model we introduce two possible solutions for the construction of $\HCd{,\surf}$-conforming discretisations. Finally, we consider two examples computed with NGSolve\footnote{www.ngsolve.org} \cite{Sch2014,Sch1997}, and discuss conclusions and outlook. 

\section{The planar relaxed micromorphic model}
Let $\con{\cdot}{\cdot}$ define scalar products and $\norm{\cdot}$ corresponding standard norms, the relaxed micromorphic model as defined in \cite{Neff2014} is given by the energy functional
\begin{align}
    I(\ud, \Pm) = \dfrac{1}{2} \int_\body& \con{\Ce\sym(\D \ud - \Pm)}{\sym(\D \ud - \Pm)} + \con{\Cc\skw(\D \ud - \Pm)}{\skw(\D \ud - \Pm)} + \con{\mathbb{C}_\mathrm{micro} \sym\Pm}{\sym\Pm} \notag \\ &  + \mu_{\mathrm{macro}} \Lc^2 (\alpha_1\norm{\sym \Curl \Pm}^2 + \alpha_2 \norm{\skw \Curl \Pm}^2) \, \dd \body - \int_\body \con{\ud}{\vb{f}} + \con{\Pm}{\bm{M}} \, \dd \body \, ,
\end{align}
where $\body \subset \R^3$, $\ud: \body \to \R^3$ is the displacement field, $\Pm : \body \to \R^{3 \times 3}$ is the microdistortion, $\vb{f}: \body \to \R^3$ are the body forces and $\bm{M}:\body \to \R^{3 \times 3}$ are couple forces. The elastic material parameters are given by the positive-definite fourth order tensors $\Ce : \R^{3 \times 3} \to \R^{3 \times 3}$ and $\mathbb{C}_\mathrm{micro}: \R^{3 \times 3} \to \R^{3 \times 3}$. The rotational coupling is given by the semi-positive-definite Cosserat \cite{boon2024mixed,Ghiba} coupling tensor $\Cc: \so(3) \to \so(3)$, $\mu_\mathrm{macro} \in \R$ is a macroscopic shear modulus, $\Lc \in \R$ is the characteristic length-scale parameter, and $\alpha_1$ and $\alpha_2$ are dimensionless weight-parameters. 
\textbf{Hereinafter we employ the abbreviations $\mathrm{M} \equiv \mathrm{macro}$ and $\mathrm{m} \equiv \mathrm{micro}$ for better readability}.
Further, for the remainder of this work we assume isotropic material behaviour 
\begin{align}
        &\Ce \bm{S}  = 2\mue \bm{S} + \lame (\tr \bm{S}) \one \, ,
        && \Cmic \bm{S} = 2\mumi \bm{S} + \lammi (\tr \bm{S}) \one \, , && \Cc \bm{A} = 2\muc \bm{A} \, , && \bm{S} = \bm{S}^T \, , && \bm{A} = -\bm{A}^T \, ,
        \label{eq:iso}
\end{align}
with $\one = \vb{e}_1 \otimes \vb{e}_1 + \vb{e}_2 \otimes \vb{e}_2 + \vb{e}_3 \otimes \vb{e}_3$, and material parameters $\muma,\mumi,\Lc>0$, $\lammi,\lame,\muc\geq 0$ such that $\Ce$ and $\Cmic$ are positive definite and $\Cc$ is positive semi definite.
Setting $\alpha_1 = 1$ and $\alpha_2 = 0$, we retrieve the $\sym \Curl$-version of the relaxed micromorphic model as defined in \cite{Lewintan2021}
\begin{align}
    I(\ud, \Pm) = \dfrac{1}{2} \int_\body& \con{\Ce\sym(\D \ud - \Pm)}{\sym(\D \ud - \Pm)} + \con{\Cc\skw(\D \ud - \Pm)}{\skw(\D \ud - \Pm)}  \notag \\ & + \con{\Cmic \sym\Pm}{\sym\Pm}  + \muma \Lc^2 \norm{\sym \Curl \Pm}^2 \, \dd \body - \int_\body \con{\ud}{\vb{f}} + \con{\Pm}{\bm{M}} \, \dd \body \, .
\end{align}
Its variation of the energy functional with respect to $\ud$ and $\Pm$ yields
\begin{align}
    &\int_\body \con{\Ce\sym\D \delta\ud}{\sym(\D \ud - \Pm)} + \con{\Cc\skw\D \delta\ud}{\skw(\D \ud - \Pm)} \, \dd \body = \int_\body \con{\delta\ud}{\vb{f}} \, \dd \body \, ,
\end{align}
and respectively
\begin{align}
    \int_\body &-\con{\Ce\sym \delta\Pm}{\sym(\D \ud - \Pm)} - \con{\Cc\skw  \delta\Pm}{\skw(\D \ud - \Pm)} + \con{\Cmic\sym\delta\Pm}{\sym\Pm} \notag \\  &+ \muma \Lc^2 \con{\sym \Curl \delta \Pm}{\sym \Curl \Pm} \, \dd \body = \int_\body  \con{\delta\Pm}{\bm{M}} \, \dd \body \, ,
\end{align}
which can be combined into
\begin{align}
    &\int_\body \con{\Ce\sym(\D \delta\ud - \delta\Pm)}{\sym(\D \ud - \Pm)} + \con{\Cc\skw(\D \delta\ud - \delta\Pm)}{\skw(\D \ud - \Pm)} + \con{\Cmic\sym\delta\Pm}{\sym\Pm} \notag \\ &\quad + \muma \Lc^2 \con{\sym \Curl \delta \Pm}{\sym \Curl \Pm} \, \dd \body = \int_\body \con{\delta\ud}{\vb{f}} + \con{\delta\Pm}{\bm{M}} \, \dd \body \, .
\end{align}
The latter is well-posed for $\{\ud, \Pm\} \in [\Honez(\body)]^3 \times \HsC{,\body}$ \cite{Lewintan2021}, where the zero index on a Hilbert space implies a vanishing Sobolev trace \cite{Hiptmair} on the boundary $\Honez(\body) = \{ u \in \Hone(\body) \; | \;  u |_{\partial \body} = 0  \}$. In \cite{SKYNOVEL,SkyOn} the characteristics of the $\HsC{,\body}$-space are explored  and it is observed that spherical fields 
\begin{align}
    &\sphP \in \Le(\body) \cdot \one \, , && \one = \vb{e}_1 \otimes \vb{e}_1 + \vb{e}_2 \otimes \vb{e}_2 + \vb{e}_3 \otimes \vb{e}_3 \in \R^{3 \times 3} \, ,
\end{align}
are in the kernel of both the $\sym\Curl$-operator and the Sobolev trace of the space 
\begin{align}
     &\forall \, \sphP \in \Le(\body) \cdot \one: && \sphP \in \ker(\sym\Curl) \cap \HsC{,\body} \, ,  && \sphP \in \ker(\tr_{\HsC{}}) \cap \HsC{,\body}   \, .
\end{align}
The direct consequence of this characteristic is that the spherical part of the microdistortion tensor $\sph \Pm = (1/3)(\tr\Pm)\one$ is allowed to be discontinuous in the domain $\body$. In a linear elastic theory, this corresponds to a jumping dilatation field, which we motivate by a complex micro-structured material. Further, we observe that for $\lammi \to \infty$ we have $\tr \Pm \to 0$ for finite energies and the microdistortion becomes deviatoric.
Unfortunately, the discontinuous dilatation property is not automatically inherited for the dimensionally reduced model of plane strain. This is because in two-dimensions $\surf \subset \R^2$ the microdistortion field is modelled as a planar field $\Pm : \surf \to \R^{2 \times 2}$ and the $\sym\Curl$-operator degenerates accordingly
\begin{align}
    & \left \| \sym \Curl \begin{bmatrix}
        \Pm & \vb{0} \\
        \vb{0} & 0 
    \end{bmatrix} \right \|^2 = \norm{\Curl \Pm}^2 = \norm{\Di (\Pm \bm{R}^T)}^2 \, , && \bm{R} = \vb{e}_1 \otimes \vb{e}_2 - \vb{e}_2 \otimes \vb{e}_1 \, ,
\end{align}
such that in two-dimensions the operator yields a vector field $\Curl \Pm : \surf \to \R^2$, see \cref{ap:b}. Consequently, the planar relaxed micromorphic model is well-posed in $\{\ud, \Pm\} \in [\Honez(\surf)]^2 \times \HC{,\surf}$. Unlike for $\HsC{,\body}$, two-dimensional spherical fields 
\begin{align}
    &\sphP \in \Le(\surf) \cdot \one \, , && \one = \vb{e}_1 \otimes \vb{e}_1 + \vb{e}_2 \otimes \vb{e}_2  \in \R^{2 \times 2} \, ,
\end{align}
are not in the kernel of the Sobolev trace operator $\ker(\tr_{\HC{}})$ corresponding with the Hilbert space $\HC{,\surf}$, and $\sph \Pm = (1/2)(\tr\Pm)\one$ is generally not allowed to jump in the domain, see \cref{ap:a}.

\begin{figure}
		\centering
		\definecolor{asl}{rgb}{0.4980392156862745,0.,1.}
		\definecolor{asb}{rgb}{0.,0.4,0.6}
		\begin{tikzpicture}[line cap=round,line join=round,>=triangle 45,x=1.0cm,y=1.0cm]
			\clip(-0.5,-0.5) rectangle (16,4.5);
			
			\fill [asb, opacity=0.1] plot [smooth cycle] coordinates {(1,3) (3,4) (7, 2) (10,3) (12,1) (10,0) (5,1) (2,1)};
			
			\begin{scope}
				\clip(5,-0.5) rectangle (12.5,1.5);
				\draw [asl, dashed] plot [smooth cycle] coordinates {(1,3) (3,4) (7, 2) (10,3) (12,1) (10,0) (5,1) (2,1)};
			\end{scope}
		    \begin{scope}
		    	\clip(5,1.5) rectangle (12.5,4.5);
		    	\draw [asl, dashed] plot [smooth cycle] coordinates {(1,3) (3,4) (7, 2) (10,3) (12,1) (10,0) (5,1) (2,1)};
		    \end{scope}
			\begin{scope}
				\clip(0,-0.5) rectangle (5,4.5);
				\draw [asb] plot [smooth cycle] coordinates {(1,3) (3,4) (7, 2) (10,3) (12,1) (10,0) (5,1) (2,1)};
			\end{scope}
			
			\draw [-to,color=black,line width=1.pt] (0,0) -- (1,0);
			\draw [-to,color=black,line width=1.pt] (0,0) -- (0,1);
			\draw (1,0) node[color=black,anchor=west] {$x$};
			\draw (0,1) node[color=black,anchor=south] {$y$};
			
			\draw [-to,color=asb,line width=1.pt] (10.,1.2) -- (10.3,0.9);
			\draw [-to,color=asb,line width=1.pt] (9.6,1.2) -- (9.3,0.9);
			\draw [-to,color=asb,line width=1.pt] (9.8,1.2) -- (9.8,0.7);
			
			\draw [-to,color=asl,line width=1.pt] (11.35,2) -- (12,2.55);
			\draw (11.8,2.3) node[color=asl,anchor=north] {$\vb{n}$};
			\draw (9.8,1.2) node[color=asb,anchor=south] {$\vb{f}$};
			\draw (2.8,2.3) node[color=asb,anchor=south] {$\bm{M}$};
			\draw [-to,asb,domain=0:180,line width=1.pt] plot ({0.5*cos(\x-180)+2.8}, {0.5*sin(\x-180)+2.5});
			
			\draw (6.5,1.15) node[color=asb,anchor=south] {$\body$};
			\draw (4.3,3.5) node[color=asb,anchor=west] {$\surf_D=\surf^u_D=\surf^P_D$};
			\draw (10.6,2.95) node[color=asl,anchor=south] {$\surf_N =\surf_N^u= \surf_N^P$};
			
			\draw [black,domain=0:360,densely dotted] plot ({0.3*cos(\x)+11.3}, {0.3*sin(\x)+1});
			\draw [black,domain=0:360,densely dotted] plot ({0.78*cos(\x)+14}, {0.78*sin(\x)+1});
			\draw [color=black, densely dotted] (11.3,1.3) -- (14,1.78);
			\draw [color=black, densely dotted] (11.3,0.7) -- (14,0.22);
			
			\draw [color=black, line width=1] (13.9,1.3) -- (14.1,1.3) -- (14.1,1.1) -- (14.3,1.1) -- (14.3,0.9) -- (14.1,0.9) -- (14.1,0.7) -- (13.9,0.7) -- (13.9,0.9)
			-- (13.7,0.9) -- (13.7,1.1) -- (13.7,1.1) -- (13.9,1.1) -- (13.9,1.3);
			\draw [color=black, line width=1] (13.7,1.7) -- (13.9,1.7) -- (13.9,1.5) -- (14.1,1.5) -- (14.1,1.7) -- (14.3,1.7);
			\draw [color=black, line width=1] (14.7,1.3) -- (14.7,1.1) -- (14.5,1.1) -- (14.5,0.9) -- (14.7,0.9) -- (14.7,0.7);
			\draw [color=black, line width=1] (13.3,1.3) -- (13.3,1.1) -- (13.5,1.1) -- (13.5,0.9) -- (13.3,0.9) --  (13.3,0.7);
			\draw [color=black, line width=1] (13.7,0.3) -- (13.9,0.3) -- (13.9,0.5) -- (14.1,0.5) -- (14.1,0.3) -- (14.3,0.3);
			
			\fill [asb,domain=0:360, opacity=0.1] plot ({0.78*cos(\x)+14}, {0.78*sin(\x)+1});
			
			\fill[white] (13.9,1.7) -- (13.9,1.5) -- (14.1,1.5) -- (14.1,1.7) -- cycle;
			
			\fill[white] (14.7,1.1) -- (14.5,1.1) -- (14.5,0.9) -- (14.7,0.9)  -- cycle;
			
			\fill[white] (13.3,1.1) -- (13.5,1.1) -- (13.5,0.9) -- (13.3,0.9) -- cycle;
			
			\fill[white] (13.9,0.3) -- (13.9,0.5) -- (14.1,0.5) -- (14.1,0.3) -- cycle;
			
			\fill[white] (13.9,1.3) -- (14.1,1.3) -- (14.1,1.1) -- (14.3,1.1) -- (14.3,0.9) -- (14.1,0.9) -- (14.1,0.7) -- (13.9,0.7) -- (13.9,0.9)
			-- (13.7,0.9) -- (13.7,1.1) -- (13.7,1.1) -- (13.9,1.1) -- cycle;
			
			\begin{scope}
				\clip(13.3,1.3) rectangle (13.1,0.7);
				\fill [white,domain=0:360] plot ({0.78*cos(\x)+14}, {0.78*sin(\x)+1});
			\end{scope}
		    \begin{scope}
		    	\clip(13.7,1.7) rectangle (14.3,2.3);
		    	\fill [white,domain=0:360] plot ({0.78*cos(\x)+14}, {0.78*sin(\x)+1});
		    \end{scope}
	        \begin{scope}
	        	\clip(13.7,0.3) rectangle (14.3,-0.3);
	        	\fill [white,domain=0:360] plot ({0.78*cos(\x)+14}, {0.78*sin(\x)+1});
	        \end{scope}
            \begin{scope}
            	\clip(14.7,1.3) rectangle (14.9,0.7);
            	\fill [white,domain=0:360] plot ({0.78*cos(\x)+14}, {0.78*sin(\x)+1});
            \end{scope}
		\end{tikzpicture}
		\caption{The domain $\body \subset \R^3$ in the relaxed micromorphic model with Neumann $\surf_N$ and Dirichlet $\surf_D$ boundaries under internal forces $\vb{f}$ and couple-forces $\bm{M}$. The Dirichlet boundary of the microdistortion is the consistent coupling condition. The model captures the averaged complex kinematics of an underlying micro-structure via the microdistortion.}
		\label{fig:domainrmm}
\end{figure}
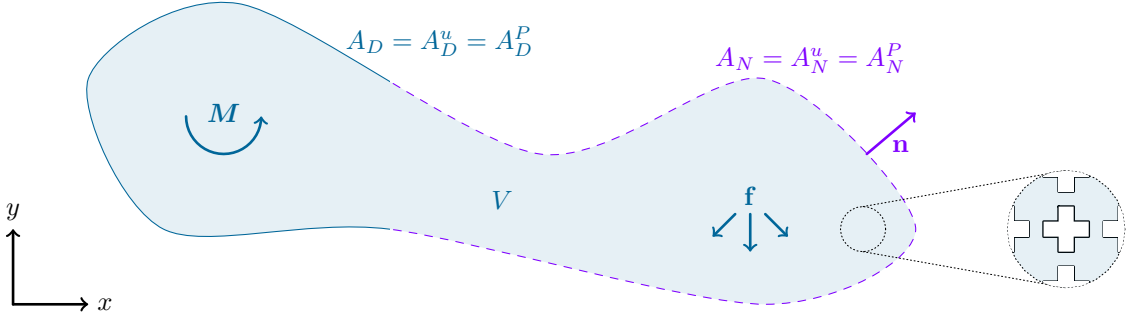

In order to restore the capacity for jumping dilatation fields also in two dimensions, we first observe that for the three-dimensional $\sym\Curl$-operator there holds
\begin{align}
    &\sym\Curl \Pm = \sym\Curl (\dev \Pm + \sph \Pm) = \sym \Curl \dev \Pm \, ,  
\end{align}
and that its corresponding Sobolev trace operator satisfies 
\begin{align}
    \tr_{\HsC{}} \Pm = \tr_{\HsC{}} (\dev \Pm + \sph \Pm) = \tr_{\HsC{}} \dev \Pm \, ,
\end{align}
since spherical fields are in the kernel of both operators. Consequently, we can decompose the $\HsC{,\body}$-space as
\begin{align}
    \HsC{,\body} = \HsCd{,\body} \oplus [\Le(\body) \cdot \one] \, ,
\end{align}
which is an orthogonal decomposition \cite{SKYNOVEL}.
This gives the option of rewriting the microdistortion field as
\begin{align}
    &\bm{P} = \devP + \sphP \, , && \devP \in \HsCd{,\body} \, , && \sphP \in \Le(\body) \cdot \one \, , 
\end{align}
without actually changing the model due to the definitions $\devP = \dev \bm{P} = \Pm - \sph \Pm$ and $\sphP = \sph \Pm$ and their orthogonality to each other $\con{\devP}{\sphP} = 0$. The variational problem now reads
\begin{align}
    &\int_\body \con{\Ce\sym(\D \delta\ud - \delta\devP - \delta\sphP)}{\sym(\D \ud - \devP - \sphP)} + \con{\Cc\skw(\D \delta\ud - \delta\devP)}{\skw(\D \ud - \devP)} \notag \\ & \qquad + \con{\Cmic\sym(\delta\devP + \delta\sphP)}{\sym(\devP + \sphP)} + \muma \Lc^2 \con{\sym \Curl \delta \devP}{\sym \Curl \devP} \, \dd \body \notag \\ & \qquad \qquad = \int_\body \con{\delta\ud}{\vb{f}} + \con{\delta\devP + \delta\sphP}{\bm{M}} \, \dd \body \, ,
\end{align}
where $\{\ud,\devP,\sphP\} \in [\Honez(\body)]^3 \times \HsCd{,\body} \times [\Le(\body) \cdot \one]$. The well-posedness of the formulation follows directly from the non-decomposed version since $\HsC{,\body} = \HsCd{,\body} \oplus [\Le(\body) \cdot \one]$.
We note that now $\devP \in \HsCd{,\body}$ can be treated in the context of the $\di \Di$-complex \cite{CRMECA_2023__351_S1_A8_0,PaulyDiv,dipietro2023discrete} (see \cref{eq:divDiv}) and make use of its corresponding finite element spaces. The boundary value problem without tractions or couple-tractions is retrieved by partial integration and splitting the boundary between Neumann and Dirichlet $\partial \body = \surf_N \cup \surf_D$
\begin{align}
    &\begin{aligned}
			    -\Di[\Ce \sym (\D \vb{u} - \devP - \sphP) + \Cc \skw (\D \vb{u} - \devP)] &= \vb{f} && \text{in} && \body \, , 
     \\
			-\dev[\Ce  \sym (\D \vb{u} - \devP) + \Cc  \skw(\D \vb{u} - \devP) - \Cmic \sym \devP - \muma \, \Lc ^ 2  \Curl (\sym \Curl\devP)] &= \dev\bm{M} && \text{in} && \body \, , 
    \\
			-\sph[\Ce  \sym (\D \vb{u} - \sphP) - \Cmic \sphP] &= \sph\bm{M} && \text{in} && \body \, ,   
			\end{aligned} \notag \\
   &\qquad\qquad\qquad \begin{aligned}
       [\Ce \sym (\D \vb{u} - \devP - \sphP) + \Cc \skw (\D \vb{u} - \devP)] \vb{n} &= 0 && \text{on} && \surf_N \,, 
    \\
		\muma \, \Lc ^ 2  (\sym \Curl\devP) (\Anti\vb{n}) &= 0 && \text{on} && \surf_N \, ,
   \\
		\vb{u} &= \widetilde{\vb{u}} && \text{on} && \surf_D \, ,
   \\
		\sym[\devP(\Anti \vb{n})] &= \sym[\D\widetilde{\vb{u}} (\Anti \vb{n})] && \text{on} && \surf_D \, ,
   \end{aligned}
   \label{eq:str3d}
\end{align}
where we employed the consistent coupling condition on the Dirichlet boundary \cite{d’Agostino2022}. The domain along with boundary conditions is depicted in \cref{fig:domainrmm}. 
Finally, the re-written energy functional now reads 
\begin{align}
    I(\ud, \devP, \sphP) = \dfrac{1}{2} \int_\body& \con{\Ce\sym(\D \ud - \devP- \sphP)}{\sym(\D \ud - \devP - \sphP)} + \con{\Cc\skw(\D \ud - \devP)}{\skw(\D \ud - \devP)}  \\ & + \con{\Cmic \sym(\devP+ \sphP)}{\sym(\devP + \sphP)} + \muma \Lc^2 \norm{\sym \Curl \devP}^2 \, \dd \body - \int_\body \con{\ud}{\vb{f}} + \con{\devP + \sphP}{\bm{M}} \, \dd \body \, . \notag
\end{align}

The split of the microdistortion field $\Pm$ into its deviatoric $\devP = \dev \Pm$ and spherical $\sphP = \sph \Pm$ components, which are now handled as explicit variables, leads to the following two-dimensional energy functional of plane strain
\begin{align}
    I(\ud, \devP, \sphP) = \dfrac{1}{2} \int_\surf& \con{\Ce\sym(\D \ud - \devP- \sphP)}{\sym(\D \ud - \devP - \sphP)} + \con{\Cc\skw(\D \ud - \devP)}{\skw(\D \ud - \devP)} \label{eq:energy2d}  \\ & + \con{\Cmic \sym(\devP+ \sphP)}{\sym(\devP + \sphP)} + \muma \Lc^2 \norm{\Curl \devP}^2 \, \dd \surf - \int_\surf \con{\ud}{\vb{f}} + \con{\devP + \sphP}{\bm{M}} \, \dd \surf \, , \notag
\end{align}
where now $\vb{u}:\surf \subset \R^2 \to  \R^2$, $\devP:\surf \subset \R^2 \to \sl(2) = \spa\{ \vb{e}_1 \otimes \vb{e}_1 - \vb{e}_2 \otimes \vb{e}_2 , \vb{e}_1 \otimes \vb{e}_2 , \vb{e}_2 \otimes \vb{e}_1 \}$ and $\sphP:\surf \subset \R^2 \to [\R \cdot \one] = \spa\{\one\}$ with $\one = \vb{e}_1 \otimes \vb{e}_1 + \vb{e}_2 \otimes \vb{e}_2 \in \R^{2\times 2}$. 
Variation of the energy functional with respect to the displacement field $\ud$ yields  
\begin{align}
    &\int_\surf \con{\Ce\sym\D \delta\ud}{\sym(\D \ud - \devP - \sphP)} + \con{\Cc\skw\D \delta\ud}{\skw(\D \ud - \devP)} \, \dd \surf = \int_\surf \con{\delta\ud}{\vb{f}} \, \dd \surf \, .
\end{align}
Analogously, for the deviatoric microdistortion field $\devP$ we find
\begin{align}
    \int_\surf &-\con{\Ce\sym \delta\devP}{\sym(\D \ud - \devP)} - \con{\Cc\skw  \delta\devP}{\skw(\D \ud - \devP)} + \con{\Cmic\sym\delta\devP}{\sym\devP} \notag \\  &+ \muma \Lc^2 \con{\Curl \delta \devP}{\Curl \devP} \, \dd \surf = \int_\surf  \con{\delta\devP}{\dev\bm{M}} \, \dd \surf \, .
\end{align}
Lastly, the variation with respect to the spherical microdistortion field $\sphP$ leads to
\begin{align}
    \int_\surf &-\con{\Ce\sym \delta\sphP}{\sym(\D \ud - \sphP)} + \con{\Cmic\sym\delta\sphP}{\sym\sphP} \, \dd \surf = \int_\surf  \con{\delta\sphP}{\sph\bm{M}} \, \dd \surf \, .
\end{align}
Combining all three variational equations yields the following planar variational form  
\begin{align}
    &\int_\surf \con{\Ce\sym(\D \delta\ud - \delta\devP - \delta\sphP)}{\sym(\D \ud - \devP - \sphP)} + \con{\Cc\skw(\D \delta\ud - \delta\devP)}{\skw(\D \ud - \devP)} \notag \\ & \qquad + \con{\Cmic\sym(\delta\devP + \delta\sphP)}{\sym(\devP + \sphP)} + \muma \Lc^2 \con{\Curl \delta \devP}{\Curl \devP} \, \dd \surf \notag \\ & \qquad \qquad = \int_\surf \con{\delta\ud}{\vb{f}} + \con{\delta\devP + \delta\sphP}{\bm{M}} \, \dd \surf \, ,
    \label{eq:redrmm}
\end{align}
where $\{\ud,\devP,\sphP\} \in [\Honez(\surf)]^2 \times \HCd{,\surf} \times [\Le(\body) \cdot \one]$. We directly observe that the discontinuous dilatation field $\sphP$ is now preserved throughout the dimensional reduction and that $\lammi \to \infty$ implies $\sphP \to 0$ as in the three-dimensional model.
Partial integration while splitting the boundary between Neumann and Dirichlet $\partial \surf = \curv_N \cup \curv_D$ now leads to the corresponding boundary value problem  
\begin{align}
    &\begin{aligned}
        -\Di[\Ce \sym (\D \vb{u} - \devP - \sphP) + \Cc \skw (\D \vb{u} - \devP)] &= \vb{f} && \text{in} \quad \surf \, , 
      \\
			-\dev[\Ce  \sym (\D \vb{u} - \devP) + \Cc  \skw(\D \vb{u} - \devP) - \Cmic \sym \devP - \muma \, \Lc ^ 2  (\D \Curl\devP)\bm{R}^T] &= \dev\bm{M} && \text{in} \quad \surf \, , 
    \\
			-\sph[\Ce  \sym (\D \vb{u} - \sphP) - \Cmic  \sphP] &= \sph\bm{M} && \text{in} \quad \surf \, , 
    \end{aligned} \notag \\
    &\qquad\qquad\qquad\begin{aligned}
        [\Ce \sym (\D \vb{u} - \devP - \sphP) + \Cc \skw (\D \vb{u} - \devP)] \vb{n} &= 0 && \text{on} \quad \curv_N \, , 
    \\
		\muma \, \Lc ^ 2  \Curl\devP  &= 0 && \text{on} \quad \curv_N \, ,
   \\
		\vb{u} &= \widetilde{\vb{u}} && \text{on} \quad \curv_D \, ,
   \\
		\devP\, \vb{t} &= (\D\widetilde{\vb{u}}) \vb{t} && \text{on} \quad \curv_D \, ,
    \end{aligned}
    \label{eq:str2d}
\end{align}
where the tangent vector on the boundary is obtained via $\vb{t} = \bm{R}^T \vb{n}$. The dimensional reduction and its non-commutative characteristic are depicted in \cref{fig:reduce}.
\begin{figure}
		\centering
		\begin{tikzpicture}[line cap=round,line join=round,>=triangle 45,x=1.0cm,y=1.0cm]
				\clip(-3.5,-5.5) rectangle (12,1.5);
				\draw (2,0) node[anchor=east] {$\displaystyle\int_\body \norm{\Pm}^2 + \norm{\sym\Curl\Pm}^2 \,\dd \body$};
                \draw [-Triangle,line width=.5pt] (0,-0.5) -- (0,-3.5);
				\draw [-Triangle,line width=.5pt] (2,0) -- (6,0);
                
                \draw (6,0) node[anchor=west] {$\displaystyle\int_\body\norm{\devP + \sphP}^2 + \norm{\sym\Curl\devP}^2\,\dd\body$};
                \draw [-Triangle,line width=.5pt] (8,-0.5) -- (8,-3.5);

                \draw (2,-4) node[anchor=east] {$\displaystyle\int_\surf \norm{\Pm}^2 + \norm{\Curl\Pm}^2 \,\dd \surf$};
                \draw [-Triangle,line width=.5pt] (2,-4) -- (6,-4);
                \draw [line width=.5pt] (3.25-0.075+1,-3.85) -- (2.75-0.075+1,-4.15);
                \draw [line width=.5pt] (3.25+0.075+1,-3.85) -- (2.75+0.075+1,-4.15);

                \draw (6,-4) node[anchor=west] {$\displaystyle\int_\surf\norm{\devP + \sphP}^2 + \norm{\Curl\devP}^2\,\dd\surf$};

                \draw (4,0) node[anchor=south] {\begin{tabular}{c}
                     $\Pm = \dev\Pm + \sph \Pm$  \\
                    $\devP = \dev\Pm$ \\
                    $\sphP = \sph\Pm$
                \end{tabular}};

                \draw (4,0) node[anchor=north] {\begin{tabular}{c}
                     decomposition  
                \end{tabular}};

                \draw (4,-4) node[anchor=north] {\begin{tabular}{c}
                     $\Pm = \dev\Pm + \sph \Pm$  \\
                    $\devP = \dev\Pm$ \\
                    $\sphP = \sph\Pm$
                \end{tabular}};

                \draw (4,-4) node[anchor=south] {(contradiction)};

                \draw (0,-2) node[anchor=west] {\begin{tabular}{c}
                     plane \\
                     strain
                \end{tabular}};
                \draw (0,-2) node[anchor=east] {\begin{tabular}{c}
                     $\bm{P}:\surf \subset \R^2 \to \R^{2 \times 2}$ 
                \end{tabular}};
                \draw (8,-2) node[anchor=east] {\begin{tabular}{c}
                     plane \\
                     strain
                \end{tabular}};
                \draw (8,-2) node[anchor=west] {\begin{tabular}{c}
                     $\devP:\surf \subset \R^2 \to \sl(2)$ \\
                     $\sphP:\surf \subset \R^2 \to \Sph$
                \end{tabular}};

                \node at (4,-2) [rectangle,draw] {\begin{tabular}{c}
                     non-  \\
                     commutative
                \end{tabular}};
			\end{tikzpicture}
		\caption{Dimensional reduction of the relaxed micromorphic model to its plane strain variant. Since the complete strain measure $\D \ud - \Pm$ is reduced to its plane strain form $(\D \ud - \Pm):\surf \subset \R^2 \to \R^{2\times 2}$, the microdistortion is defined accordingly $\Pm : \surf\subset \R^2 \to \R^{2 \times 2}$. The flow-graph highlights the non-commutativity of the procedure with respect to the isochoric decomposition $\Pm = \dev\Pm + \sph \Pm$, such that the decomposition must be applied explicitly before the dimensional reduction.}
		\label{fig:reduce}
\end{figure}

\begin{theorem}[Well-posedness of the planar form] \label{th:wellposed}
    The bilinear form $a(\{\delta\ud,\delta\devP,\delta\sphP\},\{\ud,\devP,\sphP\})$ given by the left-hand side of \cref{eq:redrmm} in conjunction with the linear form $l(\delta\ud, \delta\devP, \delta\sphP)$ given by the right-hand side of \cref{eq:redrmm} yield the problem
    \begin{align}
        a(\{\delta\ud,\delta\devP,\delta\sphP\},\{\ud,\devP,\sphP\}) = l(\delta\ud, \delta\devP, \delta\sphP) \quad \forall \, \{\delta\ud,\delta\devP,\delta\sphP\} \in \X(\surf) \, ,\label{eq:primal_form}
    \end{align}
    which is well-posed for $\{\ud,\devP,\sphP\} \in \X(\surf)$, where $\X(\surf) = [\Honez(\surf)]^2 \times \HCd{,\surf} \times [\Le(\body) \cdot \one]$ if $\Cc$ is positive definite and $\X(\surf) = [\Honez(\surf)]^2 \times \HCdz{,\surf} \times [\Le(\body) \cdot \one]$ for $\Cc = 0$ with the definition $\HCdz{,\surf} = \{\devP \in \HCd{,\surf} \; | \; \tr_{\HC{}}\Pm |_{\partial \surf} = 0\}$.
    The space $\X(\surf)$ is equipped with the norm
    \begin{align}
        \norm{\{\ud,  \devP, \sphP\}}_{\X}^2 = \norm{\ud}_{\Hone}^2 + \norm{\devP}_{\HC{}}^2 + \norm{\sphP}_{\Le}^2 \, ,
    \end{align}
    and there holds the stability estimate
    \begin{align}
        \norm{\{\ud,  \devP, \sphP\}}_{\X} \leq c \, (\norm{\vb{f}}_{\Le} + \norm{\bm{M}}_{\Le}) \, , 
    \end{align}
    where $c = c  (\Ce,\Cmic,\Cc, \muma,\Lc)$.
\end{theorem}
\begin{proof}
    The proof is given by the Lax--Milgram theorem. The continuity of the linear form is obvious and therefore omitted. Let $c > 0$ be a generic constant that may change from line to line, the continuity of the bilinear form follows via
    \begin{align}
        a(\{\delta\ud,\delta\devP,\delta\sphP\},\{\ud,\devP,\sphP\}) & \overset{CS}{\leq}  c \, (\norm{\D \delta\ud - \delta\devP - \delta\sphP}_{\Le} \norm{\D \ud - \devP - \sphP}_{\Le}  + \norm{\delta\devP + \delta\sphP}_{\Le} \norm{\devP + \sphP}_{\Le} \notag \\ 
        &\quad  + \norm{\Curl \delta \devP}_{\Le} \norm{\Curl \devP}_{\Le}) \notag \\
        & \overset{T}{\leq} c \, [(\norm{\D \delta\ud}_{\Le} +\norm{ \delta\devP}_{\Le} + \norm{\delta\sphP}_{\Le}) (\norm{\D \ud}_{\Le} + \norm{\devP}_{\Le} +\norm{ \sphP}_{\Le})   \notag \\ 
        &\quad +(\norm{\delta\devP}_{\Le} + \norm{\delta\sphP}_{\Le}) (\norm{\devP}_{\Le} + \norm{\sphP}_{\Le})  + \norm{\Curl \delta \devP}_{\Le} \norm{\Curl \devP}_{\Le}] \notag \\
        & \overset{*}{\leq} c \, [(\norm{\ud}_{\Hone} +\norm{ \delta\devP}_{\HC{}} + \norm{\delta\sphP}_{\Le}) (\norm{\ud}_{\Hone} + \norm{\devP}_{\HC{}} +\norm{ \sphP}_{\Le})   \notag \\ 
        &\quad +(\norm{\delta\devP}_{\HC{}} + \norm{\delta\sphP}_{\Le}) (\norm{\devP}_{\HC{}} + \norm{\sphP}_{\Le})  + \norm{\delta \devP}_{\HC{}} \norm{\devP}_{\HC{}}] \notag \\
        & \overset{\phantom{T}}{\leq} c \, \norm{\{\delta \ud, \delta \devP,\delta \sphP\}}_{\X} \norm{\{\ud,  \devP, \sphP\}}_{\X} \, , 
    \end{align}
    which follows from the boundedness of the fourth-order material tensors, the Cauchy-Schwarz inequality (CS), the triangle inequality (T), and the completion to the full norms (*).
    For the coercivity of the bilinear form we first assume that $\Cc$ is only positive definite. We employ again a generic constant $c > 0$ and observe
    \begin{align}
        a(\{\ud,\devP,\sphP\},\{\ud,\devP,\sphP\}) &\overset{\phantom{T}}{\geq} c \, (\norm{\D \ud - \devP - \sphP}_{\Le}^2 + \norm{\devP +\sphP}_{\Le}^2 + \norm{\Curl  \devP}_{\Le}^2)  \notag \\ 
        &\overset{\phantom{T}}{=} c \, (\norm{\D \ud}_{\Le}^2 - 2\con{\D \ud}{\devP+\sphP}_{\Le}   +\norm{\devP}^2_{\Le} +\norm{\sphP}^2_{\Le} + \norm{\devP}^2_{\Le} +\norm{\sphP}_{\Le}^2 + \norm{\Curl  \devP}_{\Le}^2)  \notag \\
        &\overset{Y}{\geq} c \, [(1-\epsilon)\norm{\D \ud}_{\Le}^2 -\epsilon^{-1}\norm{\devP+\sphP}_{\Le}^2  + 2\norm{\devP}^2_{\Le} + 2\norm{\sphP}^2_{\Le} + \norm{\Curl  \devP}_{\Le}^2]  \notag \\
        &\overset{\phantom{T}}{=} c \, [(1-\epsilon)\norm{\D \ud}_{\Le}^2  + (2-\epsilon^{-1})\norm{\devP}^2_{\Le} + (2-\epsilon^{-1})\norm{\sphP}^2_{\Le} + \norm{\Curl  \devP}_{\Le}^2]  \notag \\
        &\overset{PF}{\geq} c \, (\norm{\ud}_{\Hone}^2 + \norm{\devP}^2_{\Le} + \norm{\sphP}^2_{\Le} + \norm{\Curl  \devP}_{\Le}^2)  \notag \\
        &\overset{\phantom{T}}{\geq} c \, \norm{\{\ud,  \devP, \sphP\}}_{\X}^2  \,,
    \end{align}
    which follows by the positive-definiteness of the material tensors, the orthogonality of $\con{\devP}{\sphP}_{\Le} = 0$, Young's inequality (Y), and the Poincar\'e--Friedrich inequality (PF). The inequality is satisfied for $1/2<\epsilon <1$.
    If $\Cc$ is positive semi-definite, then we lose some control over the skew-symmetric tensors $\skw \D \ud$ and $\skw \Pm$ for $\Cc = 0$. However, for a non-vanishing Dirichlet boundary $|\curv_D| > 0$ in the variable $\devP \in \HCdz{,\surf}$, the problem is still well-posed due to 
    \begin{align}
        a(\{\ud,\devP,\sphP\},\{\ud,\devP,\sphP\}) &\overset{\phantom{T}}{\geq} c \, (\norm{ \sym(\D \ud - \devP - \sphP)}_{\Le}^2 + \norm{\sym(\devP +\sphP)}_{\Le}^2 + \norm{\Curl  \devP}_{\Le}^2)  \notag \\ 
        &\overset{\phantom{T}}{=} c \, (\norm{\sym\D \ud}_{\Le}^2 - 2\con{\sym\D \ud}{\sym\devP+\sphP}_{\Le}  + 2\norm{\sym\devP}^2_{\Le} + 2\norm{\sphP}^2_{\Le} + \norm{\Curl  \devP}_{\Le}^2)  \notag \\
        &\overset{Y}{\geq} c \, [(1-\epsilon)\norm{\sym\D \ud}_{\Le}^2  + (2-\epsilon^{-1})\norm{\sym\devP}^2_{\Le} + (2-\epsilon^{-1})\norm{\sphP}^2_{\Le} + \norm{\Curl  \devP}_{\Le}^2]  \notag \\
        &\overset{K}{\geq} c \, (\norm{\ud}_{\Hone}^2 + \norm{\sym\devP}^2_{\Le} + \norm{\sphP}^2_{\Le} + \norm{\Curl  \devP}_{\Le}^2)  \notag \\
        &\overset{G}{\geq} c \, (\norm{\ud}_{\Hone}^2 + \norm{\devP}^2_{\HC{}} + \norm{\sphP}^2_{\Le})  \notag \\
        &\overset{\phantom{T}}{\geq} c \, \norm{\{\ud,  \devP, \sphP\}}_{\X}^2  \,,
    \end{align}
    where we now used the classical Korn inequality (K) for $u$ and a generalised Korn inequality (G) \cite{LewintanInc,Lewintan2021,LewintanInc2,Gmeineder1,Gmeineder2,Neff2012} for the deviatoric microdistortion $\devP$.
\end{proof}
\begin{remark}[Determination of $\sphP$] \label{re:sphP}
    Provided the constitutive relations are isotropic as in \cref{eq:iso},  
    we observe that the third equation of either \cref{eq:str3d} in three dimensions or \cref{eq:str2d} in two dimensions allows to directly determine $\sphP$ as a function of $\D \ud$ and $\sph \bm{M}$
    \begin{align}
       (\Ce + \Cmic)\sphP = \sph \bm{M} + \sph(\Ce \sym \D \ud) \quad \iff \quad \sphP = (\Ce + \Cmic)^{-1}[\sph\bm{M} + \sph(\Ce \sym \D \ud)]  \, . 
    \end{align}
    Thus, in both cases the variable $\sphP$ can be eliminated, such that a two-field problem in $\{\ud,\devP\}$ is retrieved. However, taking the resulting $\sphP$ and inserting it into the first field equation of \cref{eq:str3d} or \cref{eq:str2d} 
    \begin{align}
        -\Di[\Ce \sym (\D \vb{u} - \devP - (\Ce + \Cmic)^{-1}[\sph\bm{M} + \sph(\Ce \sym \D \ud)]) + \Cc \skw (\D \vb{u} - \devP)] &= \vb{f} \, ,
    \end{align}
    implies that derivatives of the couple-forces $\bm{M}$ are to be considered
    \begin{align}
        \Di [(\Ce + \Cmic)^{-1}\sph\bm{M}] \, .
    \end{align}
    Consequently, we explicitly retain $\sphP$ in the equations, thus circumventing this complication. 
\end{remark}
\begin{remark}[An alternative approach via the $\di\Di$-sequence]
    An alternative approach to the explicit decomposition is by generalising the definition of the differential operator acting on the microdistortion $\Pm$. In the three-dimensional $\di\Di$-sequence \cite{PaulyDiv,CRMECA_2023__351_S1_A8_0} 
    \begin{align}
        [\Hone(\body)]^3 &\xrightarrow[]{\dev\D} \HsC{,\body} \xrightarrow[]{\sym\Curl} \HdD{,\body} \xrightarrow[]{\di\Di} \Le(\body) \, ,
        \label{eq:divDiv}
    \end{align}
    and its completion to the relaxed micromorphic sequence \cite{SKYNOVEL,SkyOn} via 
    \begin{align}
        &\ker(\sym\Curl) = \dev \D [\Hone(\body)]^3 \oplus [\Le(\body) \cdot \one] \, , && \HsC{,\body} = \HsCd{,\body} \oplus [\Le(\body) \cdot \one]  \, ,
    \end{align}
    that operator is naturally $\sym\Curl$, such that $\Pm \in \HsC{,\body}$ and $\sph \Pm$ is allowed to be discontinuous within the domain $\body \subset \R^3$. The $\di\Di$-sequence in two dimensions reads
    \begin{align}
        [\Hone(\surf)]^2 &\xrightarrow[]{\dev\D} \H^{\mathrm{dev}}(\di \Curl,\surf) \xrightarrow[]{\di \Curl} \Le(\surf) \, ,
    \end{align}
    such that one may choose $\Pm \in \H(\di \Curl,\surf) \equiv \H(\curl \Di,\surf)$. We observe that since
    \begin{align}
        \di\Curl[(\tr \Pm) \one] = \di \Di[(\tr \Pm) \one \bm{R}^T] = \di \bm{R}^T \nabla (\tr \Pm) = -\curl \nabla (\tr \Pm) = 0 \, ,
    \end{align}
    the spherical part of the microdistortion $\sph \Pm \in [\Le(\surf) \cdot \one] \subset \ker(\di \Curl) \cap \H(\di \Curl,\surf)$ is allowed to jump. However, this approach changes the regularity of the microdistortion $\Pm$ and turns the relaxed micromorphic model into a fourth order problem. Thus, we do not consider this approach further.  
\end{remark}

\subsection{The limits of the characteristic length-scale parameter $\Lc$} \label{sec:lclimits}
The relaxed micromorphic model is a true two-scale model, where scale transition is governed by the length-scale parameter $\Lc \in [0,+\infty)$. The case $\Lc = 0$ represents a highly homogeneous Cauchy continuum with classical kinematics and no couple stresses $\bm{M} = 0$. We observe that the decomposition in \cref{eq:str2d} satisfies these assumptions. Since $\bm{M} = 0$ implies 
\begin{align}
    \dev\Cc  \skw(\D \vb{u} - \devP) = \Cc  \skw(\D \vb{u} - \devP) = 0 \, ,
\end{align}
as skew-symmetric fields are naturally deviatoric, the second field equation reads
\begin{align}
    -\dev[\Ce  \sym (\D \vb{u} - \devP) - \Cmic \sym \devP] &= 0 \, ,
\end{align}
and the third field equations turns into 
\begin{align}
    -\sph[\Ce  \sym (\D \vb{u} - \sphP) - \Cmic  \sphP] &= 0 \,.  
\end{align}
The two can now be combined to find
\begin{align}
    \sym\devP + \sphP = (\Ce + \Cmic)^{-1}\Ce\sym \D \ud  \,,
\end{align}
where we used that isotropic material tensors maintain deviatoric and spherical tensors, and that the $\sym$-operator does not change tracelessness. The algebraic expression can now be inserted into the first field equation in \cref{eq:str2d} to find   
\begin{align}
    &-\Di(\Cmac\sym \D \ud) = \vb{f} \, , && \Cmac = \Cmic(\Ce + \Cmic)^{-1}\Ce \, ,
    \label{eq:lctozero}
\end{align}
which agrees with the homogenisation formula in \cite{Barbagallo2017,GOURGIOTIS2024112700} for the full Curl model.

The second limit case $\Lc \to + \infty$ represents a zoom into the micro-structure, where the micro-stiffness $\Cmic$ dominates the behaviour and no body forces occur $\vb{f} = 0$. Unlike in the previous limit, this case is more complicated and does not agree with the result of the full Curl model, compare with \cite{SKY2022115298}. Let $\Lc \to +\infty$. Then there must hold on contractible domains 
\begin{align}
    \Curl \devP = 0 \quad \iff \quad \devP = \D \pdevP \, ,
\end{align}
for finite energies in \cref{eq:energy2d}, where $\pdevP : \surf \to \R^2$ is some vectorial potential. Let the stress tensor read 
\begin{align}
    \bm{\sigma} &= 2\mue \sym(\D \ud  - \Pm) + 2\muc \skw(\D \ud  - \Pm) + \lame \tr(\D \ud  - \Pm) \one \, , 
\end{align}
then the first field equation from \cref{eq:str2d} without body forces can be written as
\begin{align}
    -\Di \bm{\sigma} = -\Di[\dev\Ce \sym (\D \vb{u} - \devP) + \sph\Ce \sym (\D \vb{u} - \sphP) + \Cc \skw (\D \vb{u} - \devP)] &= 0 \, ,\label{eq:field_a_lc}
\end{align}
where we applied the isochoric decomposition on the the stress tensor $\bm{\sigma} = \dev\bm{\sigma} + \sph\bm{\sigma}$. Now, taking the divergence of the second and third field equations  in \cref{eq:str2d}, and subsequently combining them yields
\begin{align}
    -\Di[\dev\Ce  \sym (\D \vb{u} - \devP) + \sph\Ce  \sym (\D \vb{u} - \sphP) + \Cc  \skw(\D \vb{u} - \devP) - \Cmic (\sym \devP + \sphP)]&= \Di\bm{M}  \, ,  
\end{align}
where we recognise the first three terms to vanish by \cref{eq:field_a_lc}, resulting in 
\begin{align}
    \Di[\Cmic (\sym \devP + \sphP)]&= \Di\bm{M}  \, . \label{eq:red_eq_lc} 
\end{align}
Contrary to the full Curl formulation of the relaxed micromorphic model that yields the micro-level elasticity problem $\Di(\Cmic \sym \Pm) = \Di(\Cmic \sym \D \vb{p}) = \Di\bm{M}$ with some potential vector $\vb{p}:\surf \to \R^2: \Curl \Pm = \Curl \D \vb{p} = 0$, we find here a two-field problem. However, since $\devP$ is traceless, its vectorial potential must satisfy 
\begin{align}
    \tr \devP = \tr \D \pdevP = \di \pdevP = 0 \, ,
\end{align}
which is exactly satisfied by a two-dimensional scalar curl potential $\mathfrak{p} : \surf \to \R$   
\begin{align}
    &\pdevP = \bm{R} \nabla \mathfrak{p}  \, , && \di \bm{R} \nabla \mathfrak{p}= 0  \, . \label{eq:kernel}   
\end{align}
Further, we can rewrite $\sphP$ as
\begin{align}
    &\sphP = \dfrac{1}{2} \, \mathfrak{j} \, \one \, , && \mathfrak{j}:\surf \to \R \, ,
\end{align}
where $\one$ is the constant identity tensor and $\mathfrak{j}$ is the underlying function. Put together along with the explicit constitutive equation, the two-field problem \cref{eq:red_eq_lc} can be recast as
\begin{align}
    \Di[ 2\mumi\sym \D\bm{R}\nabla\mathfrak{p} + (\mumi + \lammi) \, \mathfrak{j} \, \one] =  2\mumi\Delta \bm{R}\nabla\mathfrak{p} + (\mumi + \lammi) \nabla \mathfrak{j}&= \Di\bm{M}  \, . \label{eq:step}   
\end{align}
Finally, taking the curl of the latter leads to the fourth-order, two-dimensional quad-curl problem \cite{DiPietroRot,Sun2019,WangLiZhang} 
\begin{align}
       2\mumi \curl \Delta \bm{R}\nabla\mathfrak{p} = -2\mumi \curl \bm{R}\nabla \curl  \bm{R}\nabla\mathfrak{p} &= \curl \Di\bm{M} \, ,  \label{eq:quad_curl_lc}
\end{align}
where we used $\Delta(\cdot) = \di \nabla(\cdot) = -\curl \bm{R} \nabla(\cdot)$. Thus, \cref{eq:quad_curl_lc} can also be expressed as a biharmonic equation
\begin{align}
    -2\mumi \Delta \Delta \mathfrak{p} &= \curl \Di\bm{M} \, . 
\end{align}
Analogously, a classical Poisson problem is retrieved from \cref{eq:step} by taking the divergence on both sides
\begin{align}
     (\mumi + \lammi) \Delta \mathfrak{j}=  \di\Di\bm{M}  \, ,\label{eq:poisson_j}
\end{align}
where we used $2\mumi \di (\Delta \bm{R} \nabla \mathfrak{p}) = 2\mumi \Delta \di( \bm{R} \nabla \mathfrak{p}) = 2\mumi \Delta \curl \nabla \mathfrak{p} = 0$.
At this point we recall that by as per \cref{re:sphP} $\sphP$ can be algebraically expressed as a function of $\ud$ and $\bm{M}$
\begin{align}
    \sphP = (\Ce + \Cmic)^{-1}(\sph\bm{M} + \sph\Ce \sym \D \ud) = \dfrac{1}{2[(\mue + \mumi) + (\lame + \lammi)] } \left[ \dfrac{1}{2} (\tr \bm{M}) \one + (\mue + \lame) (\di \ud) \one \right  ] \, .
\end{align}
Consequently, we can identify $\mathfrak{j}$ as 
\begin{align}
    \mathfrak{j} = \tr\sphP = \dfrac{1}{2[(\mue + \mumi) + (\lame + \lammi)] } [  \tr \bm{M} + 2(\mue + \lame) \di \ud ] \, .\label{eq:j_expr}
\end{align}
Multiplying \cref{eq:j_expr} with $(\mumi + \lammi)$ yields
\begin{align}
    (\mumi + \lammi) \,\mathfrak{j} = \dfrac{\mumi + \lammi}{2[(\mue + \mumi) + (\lame + \lammi)] } [  \tr \bm{M} + 2(\mue + \lame) \di \ud ] = \dfrac{(\muma + \lamma)}{2} \tr \bm{M} + (\muma + \lamma) \di \ud   \, ,
\end{align}
as per the homogenisation formula in \cite{GOURGIOTIS2024112700}, see \cref{eq:homogen}.
Thus, the Poisson equation \cref{eq:poisson_j} for $\mathfrak{j}$ can be recast as
\begin{align}
     -(\muma + \lamma) \Delta \di \ud =  \dfrac{(\muma + \lamma)}{2} \Delta \tr \bm{M} - \di\Di\bm{M} \, ,
\end{align}
and the source-intensity of the displacement field $\di \ud$ at the micro-level can be directly determined as
\begin{align}
     \di \ud = \dfrac{\mumi + \lammi}{\muma + \lamma} \, \widetilde{\mathfrak{j}} - \dfrac{1}{2} \tr \bm{M}   \, ,
\end{align}
for a precomputed $\widetilde{\mathfrak{j}}$, retrieved from the Poisson equation \cref{eq:poisson_j}.

Variationally, we can approach the problem $\Lc \to +\infty$ with a mixed formulation. We introduce the hyper-stress variable 
\begin{align}
     \hP = \muma\Lc^2 \Curl \devP \, , 
\end{align}
such that the variational problem in \cref{th:wellposed} is adapted to
\begin{align}
    a(\{\delta \ud, \delta \devP, \delta \sphP\}, \{ \ud,  \devP,  \sphP\}) + b(\{\delta \ud, \delta \devP, \delta \sphP\},  \hP) &= l(\{\delta \ud, \delta \devP, \delta \sphP\}) && \forall \, \{\delta \ud, \delta \devP, \delta \sphP\} \in \X(\surf)  \, , \notag \\
    b(\delta \hP, \{\ud,  \devP, \sphP\}) - \dfrac{1}{\muma \Lc^2} c(\delta \hP, \hP) &= 0 && \forall \, \delta \hP \in \Z(\surf) \, ,
    \label{eq:probmixed}
\end{align}
where $\Z(\surf)$ is to be subsequently defined, and the bilinear forms read
\begin{align}
    a(\{\delta \ud, \delta \devP, \delta \sphP\}, \{ \ud,  \devP,  \sphP\}) &= \int_\surf \con{\Ce\sym(\D \delta\ud - \delta\devP - \delta\sphP)}{\sym(\D \ud - \devP - \sphP)} 
    \notag \\ &\qquad + \con{\Cmic\sym(\delta\devP + \delta\sphP)}{\sym(\devP + \sphP)} + \con{\Cc\skw(\D \delta\ud - \delta\devP)}{\skw(\D \ud - \devP)}   \, \dd \surf \,, \notag \\
    b(\{\delta \ud, \delta \devP, \delta \sphP\},  \hP) &= \int_\surf \con{\Curl \delta \devP}{\hP} \, \dd \surf \, , \notag \\
    c(\delta \hP, \hP) &= \int_\surf \con{\delta \hP}{\hP} \, \dd \surf \, .
\end{align}
Observe that \cref{eq:kernel} implies 
\begin{align}
    \HCd{,\surf} \cap \ker(\Curl) = \D \bm{R} \nabla [\H^2(\surf)] \, , 
\end{align}
and by the two-dimensional de Rham sequence \cite{arnold_complexes_2021,PaulyDeRham}
\begin{align}
    \Hone(\surf) &\xrightarrow[]{\bm{R}\nabla} \Hc{,\surf} \xrightarrow[]{\curl} \Le(\surf) \, ,
\end{align}
the Curl-operator is a row-wise rotated divergence such that the problem can be identified with the two-dimensional grad-curl \cite{Kaiboqcurl,arnold_complexes_2021} sub-sequence 
\begin{align}
    \H^2(\surf) &\xrightarrow[]{\D\bm{R}\nabla} \HCd{,\surf} \xrightarrow[]{\Curl} [\Le(\surf)]^2 \, ,
\end{align}
where the last identity is a surjection.
Consequently, the space of the hyper-stress $\hP \in \Z(\surf)$ is $\Z(\surf) = [\Le(\surf)]^2$, or $\Z(\surf) = [\Lez(\surf)]^2$ in the case of full Dirichlet boundary conditions for $\devP$
\begin{align}
    \H^2_0(\surf) &\xrightarrow[]{\D\bm{R}\nabla} \HCdz{,\surf} \xrightarrow[]{\Curl} [\Lez(\surf)]^2 \, .
\end{align}
\begin{theorem}[Well-posedness of the mixed form]
    The problem defined in \cref{eq:probmixed} is well-posed and there holds the stability estimate
    \begin{align}
        \norm{\{\ud, \devP,\sphP\}}_X + \norm{\hP}_{\Le} \leq c \, (\norm{\vb{f}}_{\Le} + \norm{\bm{M}}_{\Le}) \, , 
    \end{align}
    where $c = c(\Ce,\Cmic,\Cc, \muma)$, such that it does not depend on $\Lc$.
\end{theorem}
\begin{proof}
    The proof follows by the extended Brezzi-theorem \cite[Thm. 4.11]{Bra2013}. The continuity of the bilinear forms is obvious and the kernel-coercivity of $a(\cdot,\cdot)$ on the sub-space $\X(A) \cap \ker(b)=\{\{ \ud,  \devP,  \sphP\}\in \X(A) \; | \; b(\{ \ud,  \devP,  \sphP\},\hP)=0 \quad \forall \, \hP\in Z(A)\}=\{\{ \ud,  \devP,  \sphP\}\in \X(\surf) \; | \; \Curl \devP=0\}$ follows from \cref{th:wellposed}. Thus, we only need to satisfy the Ladyzhenskaya--Babu{\v{s}}ka--Brezzi (LBB) condition. We choose $\delta\ud = 0$ and $\delta\sphP = 0$, as well as $\delta \devP$ such that $\Curl\delta \devP = \hP$ and $\norm{\delta\devP}_{\Le}\leq c\, \norm{\hP}_{\Le}$ by the grad-curl sequence, such that there holds
    \begin{align}
        \sup_{\{\delta \ud, \delta \devP, \delta \sphP\} \in \X} \dfrac{b(\{\delta \ud, \delta \devP, \delta \sphP\},\hP)}{\norm{\{\delta \ud, \delta \devP, \delta \sphP\}}_\X} \geq \dfrac{\int_\surf \con{\Curl \delta \devP}{\hP} \, \dd \surf}{\norm{\delta\devP}_{\Le} + \norm{\Curl \delta\devP}_{\Le}} \geq c \, \dfrac{\norm{\hP}^2_{\Le}}{\norm{\hP}_{\Le}} = c \, \norm{\hP}_{\Le}\, ,
    \end{align}
    finishing the proof.
\end{proof}

\section{Finite element discretisations}
Since the primal formulation \cref{eq:primal_form} is shown to be well-posed by the Lax--Milgram theorem, any conforming subspace discretisation is also well-posed. Conforming subspaces for $\Hone(\surf)$ and $\Le(\surf)$ are the classical continuous and discontinuous Lagrange elements of polynomial degree $p$, $\CG^p(\surf) \subset \Hone(\surf)$, $\DG^p(\surf) \subset \Le(\surf)$. Consequently, we can discretise the displacement and dilatation fields as $\{\ud, \sphP\} \in [\CG^p(\surf)]^2 \times [\DG^{p-1} \cdot \one]$, where we choose the polynomial degree $p-1$ for the dilatation field in order to match with the polynomial degree of $\D \ud$. The challenge is to construct a minimally regular conforming finite element for $\HCd{,\surf}$, which for the mixed formulation, must also commute on the grad-curl complex in order to satisfy Fortin's criterion \cite{Bra2013,arnold_complexes_2021} and ensure stable discretisations. 

An $\HCd{,\surf}$-conforming element cannot be constructed on standard triangles because the deviatoric constraint ties two tangent vectors together at a vertex, such that all interfacing elements at the vertex must share the same tangent vectors. Clearly, this is impossible on general triangulations. If one abandons the minimal regularity and enforces $\C^0(\surf)$-continuity at the vertices, it becomes simple to introduce a conforming element. We propose such a construction in \cref{sec:dev} based on the polytopal approach \cite{sky_polytopal_2022}. However, the minimal regularity is an important property in the face of jumping material coefficients, compare \cite{SKYNOVEL}, which is perfectly plausible in micro-structured materials. As such, we also present a macro element construction on Clough--Tocher splits \cite{Philippe} to derive a minimally regular, fully deviatoric $\HCd{,\surf}$-conforming finite element space, see \cref{sec:macro}.   
Unfortunately, the implementation of macro elements is not trivial, leading us to present an alternative with a weakly deviatoric formulation and N\'ed\'elec elements \cite{sky_polytopal_2022,sky_higher_2023,Ned2,Nedelec1980,haubold2023high}.  

In the following we define the reference element as 
\begin{align}
    \Gamma = \{ (\xi,\eta) \in [0,1]^2 \; | \; \xi + \eta \leq 1 \} \, , 
\end{align}
and an arbitrary element in the physical mesh as $T \subset \R^2$. Tangent and normal vectors on the reference element are denoted with $\bm{\tau}$ and $\bm{\nu}$, respectively. Their counterparts on the physical element $T$ are $\vb{t}$ and $\vb{n}$. In general we define 
\begin{align}
    &\bm{\nu} = \bm{R} \bm{\tau} \, , && \vb{n} = \bm{R} \vb{t} \, , && \bm{R} = \vb{e}_1 \otimes \vb{e}_2 - \vb{e}_2 \otimes \vb{e}_1 \, , 
\end{align}
such that 
\begin{align}
     &\norm{\bm{\tau}} = \norm{\bm{\nu}} \, , && \norm{\vb{t}} = \norm{\vb{n}} \, .
\end{align}
Tangent vectors are mapped from the reference to the physical element via the Jacobi matrix $\bm{J}= \D \vb{x}$ of the mapping $\vb{x}:\Gamma \to \surf$ via
\begin{align}
    &\vb{t} = \bm{J} \bm{\tau} \, .
\end{align}

\subsection{A strongly deviatoric element} \label{sec:dev}
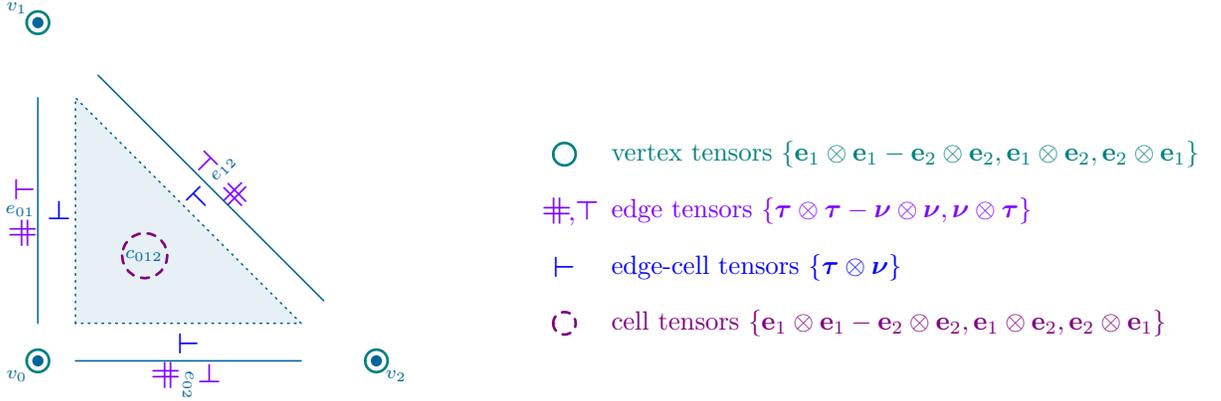
\begin{figure}
		\centering
		\definecolor{asl}{rgb}{0.4980392156862745,0.,1.}
    		\definecolor{asb}{rgb}{0.,0.4,0.6}
    		\begin{tikzpicture}[line cap=round,line join=round,>=triangle 45,x=1.0cm,y=1.0cm]
    			\clip(-2,-1) rectangle (15,4.5);
    			\draw (-0.5,-0.5) node[circle,fill=asb,inner sep=1.5pt] {};
    			\draw (-0.5,4) node[circle,fill=asb,inner sep=1.5pt] {};
    			\draw (4,-0.5) node[circle,fill=asb,inner sep=1.5pt] {};
    			\draw [color=asb,line width=.6pt] (-0.5,0) -- (-0.5,3);
    			\draw [color=asb,line width=.6pt] (0,-0.5) -- (3,-0.5);
    			\draw [color=asb,line width=.6pt] (0.3,3.3) -- (3.3,0.3);
    			\draw [dotted,color=asb,line width=.6pt] (0,0) -- (0,3) -- (3,0) -- (0,0);
    			\fill[opacity=0.1, asb] (0,0) -- (0,3) -- (3,0) -- cycle;
    			\draw (-0.5,-0.5) node[color=asb,anchor=north east] {$_{v_0}$};
    			
                \draw [teal,domain=0:360,line width=1.pt] plot ({0.15*cos(\x-180)-0.5}, {0.15*sin(\x-180)-0.5});
    			
    			\draw (4,-0.5) node[color=asb,anchor=north west] {$_{v_2}$};
    			\draw (-0.5,4) node[color=asb,anchor=south east] {$_{v_1}$};

                \draw [teal,domain=0:360,line width=1.pt] plot ({0.15*cos(\x-180)-0.5}, {0.15*sin(\x-180)+4});
                \draw [teal,domain=0:360,line width=1.pt] plot ({0.15*cos(\x-180)+4}, {0.15*sin(\x-180)-0.5});
    			
    			\draw (-0.4,1.5) node[color=asb,anchor=east] {$_{e_{01}}$};
                \draw (-0.5,1.5) node[color=blue,anchor=west] {$\bm{\perp}$};
                \draw (-0.45,1.5) node[color=asl,anchor=north east] {\rotatebox{180}{$\bm{\parallel}$}};
                \draw (-0.4,1.44) node[color=asl,anchor=north east] {\rotatebox{270}{$\bm{\parallel}$}};
                \draw (-0.43,1.5) node[color=asl,anchor=south east] {\rotatebox{-90}{$\bm{\perp}$}};
    			
    			\draw (1.5,-.5) node[color=asb,anchor=north] {\rotatebox{-90}{$_{e_{02}}$}};

                \draw (1.5,-.55) node[color=blue,anchor=south] {\rotatebox{270}{$\bm{\perp}$}};
                \draw (1.5,-.44) node[color=asl,anchor=north east] {\rotatebox{90}{$\bm{\parallel}$}};
                \draw (1.45,-.39) node[color=asl,anchor=north east] {\rotatebox{180}{$\bm{\parallel}$}};
                \draw (1.5,-.42) node[color=asl,anchor=north west] {\rotatebox{0}{$\bm{\perp}$}};
    			
    			\draw (1.66,1.74) node[color=asb,anchor=south west] {\rotatebox{45}{$_{e_{12}}$}};

                \draw (2-0.03,2-0.03) node[color=blue,anchor=north east] {\rotatebox{225}{$\bm{\perp}$}};
                \draw (1.8,1.4) node[color=asl,anchor=south west] {\rotatebox{225}{$\bm{\parallel}$}};
                \draw (1.8,1.4) node[color=asl,anchor=south west] {\rotatebox{135}{$\bm{\parallel}$}};
                \draw (1.4,1.8) node[color=asl,anchor=south west] {\rotatebox{135}{$\bm{\perp}$}};
    			
    			\draw (0.92,0.9)
    			node[color=asb] {$_{c_{012}}$};
                \draw [violet,domain=0:360,line width=1.pt, dashed] plot ({0.3*cos(\x-180)+0.92}, {0.3*sin(\x-180)+0.9});
    			
                \draw [teal,domain=0:360,line width=1.pt] plot ({0.15*cos(\x-180)+6.5}, {0.15*sin(\x-180)+2.25});
    			\draw (7,2.25)
    			node[color=teal,anchor=west] {vertex tensors $\{ \vb{e}_1 \otimes \vb{e}_1 - \vb{e}_2 \otimes \vb{e}_2, \vb{e}_1 \otimes \vb{e}_2, \vb{e}_2 \otimes \vb{e}_1\}$};
                    \draw (6.4,1.5) node[color=asl] {\rotatebox{90}{$\bm{\parallel}$}};
    			\draw (6.4,1.5) node[color=asl] {$\bm{\parallel}$};
                \draw (6.6,1.35) node[color=asl] {,};
                 \draw (6.8,1.5) node[color=asl] {\rotatebox{180}{$\bm{\perp}$}};
    			\draw (7,1.5)
    			node[color=asl,anchor=west] {edge tensors $\{ \bm{\tau} \otimes \bm{\tau} - \bm{\nu} \otimes \bm{\nu}, \bm{\nu} \otimes \bm{\tau}\}$};
                \draw (6.5,0.75) node[color=blue] {\rotatebox{270}{$\bm{\perp}$}};
    			\draw (7,0.75)
    			node[color=blue,anchor=west] {edge-cell tensors $\{\bm{\tau} \otimes \bm{\nu}\}$};
                \draw [violet,domain=0:360,line width=1.pt, dashed] plot ({0.15*cos(\x-180)+6.5}, {0.15*sin(\x-180)});
    			\draw (7,0)
    			node[color=violet,anchor=west] {cell tensors $\{ \vb{e}_1 \otimes \vb{e}_1 - \vb{e}_2 \otimes \vb{e}_2, \vb{e}_1 \otimes \vb{e}_2, \vb{e}_2 \otimes \vb{e}_1\}$};
    		\end{tikzpicture}
		\caption{Base tensors for the reference triangle element on their corresponding polytope.}
		\label{fig:template}
\end{figure}

In order to construct a strongly deviatoric $\HCd{,\surf}$-conforming element we abandon the minimal tangential regularity of $\Hc{,\surf}$ and enforce full $\C^0(\surf)$-continuity at the vertices of the triangles. For simplicity we base the construction on the barycentric coordinates 
\begin{align}
    &\lambda_0 = 1 - \xi - \eta \, , && \lambda_1 = \eta \, , && \lambda_2 = \xi \, ,
\end{align}
but note that any $\C^0(\surf)$-continuous Lagrangian space $\CG^p(\Gamma)$ can be used instead.
\begin{definition}[Strongly deviatoric base functions]
    The base functions are defined per polytope of the reference triangle.
    \begin{itemize}
        \item on each vertex $v_i$ with scalar base function $\lambda_i$ we define the base functions
        \begin{align}
            &\bm{\Upsilon}_{i j} = \lambda_i \bm{\Psi}_j \, , &&  \bm{\Psi}_j \in \{ \vb{e}_1 \otimes \vb{e}_1 - \vb{e}_2 \otimes \vb{e}_2, \vb{e}_1 \otimes \vb{e}_2, \vb{e}_2 \otimes \vb{e}_1\} \, ,
        \end{align}
        such that there are $3$ base functions for each vertex. 
        \item on each edge $e_{ij}$ where $(i,j) \in \{(0,1),(0,2),(1,2)\}$, equipped with the tangent $\bm{\tau}_{ij}$ and normal $\bm{\nu}_{ij}$ vectors we define
        \begin{align}
            &\bm{\Upsilon}_{ijk}^{qr} = \lambda_i^q \lambda_j^r  \bm{\Psi}_k \, , &&  \bm{\Psi}_k \in \{ \bm{\tau}_{ij} \otimes \bm{\tau}_{ij} - \bm{\nu}_{ij} \otimes \bm{\nu}_{ij}, \bm{\nu}_{ij} \otimes \bm{\tau}_{ij} \} \, ,
        \end{align}
        where $q,r \geq 1$. For each edge there are $2 \cdot (p-1)$ base functions, where $p$ is the polynomial order of the edge.
        \item for the cell we define the edge-cell and pure cell base functions
        \begin{align}
           &\bm{\Upsilon}_{ijk}^{qr} = \lambda_i^q \lambda_j^r  \bm{\Psi}_k \, , &&  \bm{\Psi}_k \in \{ \bm{\tau}_{ij} \otimes \bm{\nu}_{ij} \} \, , \notag \\
            &\bm{\Upsilon}_{k}^{qrz} = \lambda_0^q \lambda_1^r \lambda_2^z \bm{\Psi}_k \, , &&  \bm{\Psi}_k \in \{ \vb{e}_1 \otimes \vb{e}_1 - \vb{e}_2 \otimes \vb{e}_2, \vb{e}_1 \otimes \vb{e}_2, \vb{e}_2 \otimes \vb{e}_1\} \, ,
        \end{align} 
        where $(i,j) \in \{(0,1),(0,2),(1,2)\}$ iterate over all edges and $q,r,z \geq 1$. There are $3 \cdot (p-1)$ edge-cell base functions and $3 \cdot (p-2)(p-1)/2$ cell base functions for a uniform polynomial order $\Po^p(\Gamma)$.  
    \end{itemize}
\end{definition}
The vertex and cell base functions of the construction are defined using a global deviatoric Cartesian basis. As such, they do not undergo a transformation from the reference to the physical element. For the edge and edge-cell functions we define the  edge-wise mapping from a reference to a physical element via
\begin{align}
    & \bm{Y} = \bm{T}_{ij}\bm{\Upsilon}\bm{T}_{ij}^T \,, && \bm{T}_{ij} = \dfrac{1}{\norm{\bm{\tau}_{ij}}^2} (\vb{t}_{ij} \otimes \bm{\tau}_{ij} + \vb{n}_{ij} \otimes \bm{\nu}_{ij}) \, ,
    \label{eq:trans}
\end{align}
for each edge $e_{ij}$. The definition of base tensors on the reference triangle is illustrated in \cref{fig:template}. 

\begin{theorem}[Linear independence and conformity]
    The presented base functions are linearly independent and can be used to span an $\HCd{,\surf}$-conforming subspace.
    \label{th:strongind}
\end{theorem}
\begin{proof}
    The proof of linear independence is obvious by construction. It follows directly from the linear independence of the scalar base functions and the multiplication with linearly independent sets of constant tensors. For conformity we observe that the vertex base functions are $\C^0(\surf)$-continuous and deviatoric. On the edges connectivity is only implied for $(\cdot) \otimes \vb{t}$-type components. In other words, tangential continuity is implied, thus satisfying a vanishing jump in the Sobolev trace of $\HC{,\surf}$. Further, the edge basis is deviatoric since 
    \begin{align}
        &\tr(\vb{t} \otimes \vb{t} - \vb{n} \otimes \vb{n}) = \norm{\vb{t}}^2 - \norm{\vb{n}}^2 = 0 \, , && \tr(\vb{n} \otimes \vb{t}) = \con{\vb{n}}{\vb{t}} = 0 \, , 
    \end{align}
    as $\norm{\vb{t}} = \norm{\vb{n}}$ and $\vb{t} \perp \vb{n}$. This follows analogously for the edge-cell base functions. Finally, the pure cell base functions vanish on the entire boundary of the triangle and are built with a deviatoric basis. 
\end{proof}
Let $\devP \in [\H^{1}(\surf)\otimes\sl(2)]$ be part of the exact solution and $\Y^p(\surf)$ be the discrete space given by the former construction. Following \cite{hu_family_2014,hu_finite_2016} for the Hu--Zhang elements, we define a canonical interpolation operator $\Pi_y^p:[\H^{1}(\surf)\otimes\sl(2)]\to \Y^p(\surf)$ by slightly modifying the Scott--Zhang interpolation operator \cite{SZ1990}. We have the vertex evaluations 
\begin{align}
    &\con{\bm{V}_i}{\devP - \Pi_y^p \devP}\at_{v_j} = 0 \, , && \bm{V}_i \in \{ \vb{e}_1 \otimes \vb{e}_1 - \vb{e}_2 \otimes \vb{e}_2, \vb{e}_1 \otimes \vb{e}_2, \vb{e}_2 \otimes \vb{e}_1\} \, ,
\end{align}
if $\devP$ is continuous, or averages over edges as per \cite{SZ1990} otherwise. For the remaining polytopes we can define projections. We construct edge integrals for every edge interval $\curv_{ij}$ 
\begin{align}
    &\int_{\curv_{ij}} q_k \con{\bm{V}_l}{\devP - \Pi_y^p\devP} \, \dd \curv = 0 \, , && q_k \in \Po^{p-2}(\curv_{ij}) \, , && \bm{V}_l \in \{ \vb{t}_{ij} \otimes \vb{t}_{ij} - \vb{n}_{ij} \otimes \vb{n}_{ij}, \vb{n}_{ij} \otimes \vb{t}_{ij} \} \, ,
\end{align}
edge-cell integrals
\begin{align}
     &\int_{\curv_{ij}} q_k \con{\bm{V}_l}{\devP - \Pi_y^p\devP} \, \dd \curv = 0 \, , && q_k \in \Po^{p-2}(\curv_{ij}) \, , && \bm{V}_l \in \{\vb{t}_{ij} \otimes \vb{n}_{ij} \} \, ,
\end{align}
and finally cell integrals on every element $\elem_i$
\begin{align}
    &\int_{\elem_i} q_j \con{\bm{V}_k}{\devP - \Pi_y^p\devP}  \, \dd \surf = 0 \, , && q_j \in \Po^{p-3}(\elem_i) \, , && \bm{V}_k \in \{ \vb{e}_1 \otimes \vb{e}_1 - \vb{e}_2 \otimes \vb{e}_2, \vb{e}_1 \otimes \vb{e}_2, \vb{e}_2 \otimes \vb{e}_1\} \, .
\end{align}
The stability of $\Pi_y^p$ follows by the stability of the Scott--Zhang operator
\begin{align}
    \|\Pi_y^p\devP\|_{\HC{}}\leq C \|\devP\|_{\Hone}\,.
\end{align}
Note that the curl operator is well-defined on each element by Stokes' theorem and our degrees of freedom
\begin{align}
    \int_T \Curl (\Pi_y^p \devP) \, \dd \surf = 
\int_{\partial \elem} (\Pi_y^p \devP) \vb{t} \, \dd \curv = \sum_i \int_{\curv_i} ( \Pi_y^p \devP) \vb{t} \, \dd \curv \, .
\end{align}
Assume that $\devP \in [\H^{p+1}(\surf)\otimes\sl(2)]$. As $\Pi_y^p$ preserves locally deviatoric polynomials of order $p$ we directly obtain the approximation properties
\begin{align}
    \norm{\devP - \Pi^p_y\devP}_{\Le} \leq c\,h^{p+1}|\devP|_{\H^{p+1}} \, ,\qquad \norm{\Curl\devP - \Curl\Pi^p_y\devP}_{\Le} \leq c\,h^{p}|\devP|_{\H^{p+1}}\,.
    \label{eq:inty}
\end{align}
Consequently, analogously to the derivation in \cite{hu_family_2014,hu_finite_2016} for the elasticity sub-complex \cite{arnold_mixed_2002,arnold_complexes_2021,sky2023reissnermindlin,SKYuni}, we get that our element allows for a commuting grad-curl sub-complex for the polynomial degree $p \geq 3$
\begin{align}
    \begin{matrix}
        \H^2(\surf) &\xrightarrow[]{\D\bm{R}\nabla} &\HCd{,\surf} &\xrightarrow[]{\Curl} &[\Le(\surf)]^2 \\[0.5em]
        \Pi_a^5\bigg\downarrow & &\Pi_y^3\bigg\downarrow & &\Pi_o^2\bigg\downarrow \\[1em]
        \mathcal{A}^5(\surf) &\xrightarrow[]{\D\bm{R}\nabla} &\Y^3(\surf) &\xrightarrow[]{\Curl} &[\DG^2(\surf)]^2 
    \end{matrix} \quad \, ,
\end{align}
where $\mathcal{A}^5(\surf)$ is the quintic $\C^1(\surf)$-continuous Argyris element \cite{CarstensenHu,Dominguez,kirby_general_2018}. 

\subsection{A macro element on Clough--Tocher splits} \label{sec:macro}
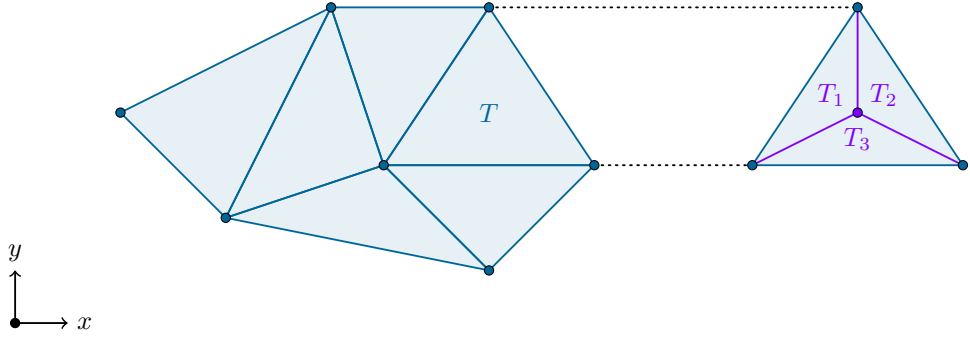
\begin{figure}
    \centering
    \definecolor{xfqqff}{rgb}{0.4980392156862745,0,1}
\definecolor{qqwwzz}{rgb}{0,0.4,0.6}
\begin{tikzpicture}[scale = 0.7, line cap=round,line join=round,>=triangle 45,x=1cm,y=1cm]
\clip(1.5,0.5) rectangle (20.5,7.5);
\fill[line width=0.7pt,color=qqwwzz,fill=qqwwzz,fill opacity=0.1] (4,5) -- (8,7) -- (6,3) -- cycle;
\fill[line width=0.7pt,color=qqwwzz,fill=qqwwzz,fill opacity=0.1] (6,3) -- (9,4) -- (8,7) -- cycle;
\fill[line width=0.7pt,color=qqwwzz,fill=qqwwzz,fill opacity=0.1] (9,4) -- (11,2) -- (6,3) -- cycle;
\fill[line width=0.7pt,color=qqwwzz,fill=qqwwzz,fill opacity=0.1] (8,7) -- (11,7) -- (9,4) -- cycle;
\fill[line width=0.7pt,color=qqwwzz,fill=qqwwzz,fill opacity=0.1] (13,4) -- (9,4) -- (11,2) -- cycle;
\fill[line width=0.7pt,color=qqwwzz,fill=qqwwzz,fill opacity=0.1] (13,4) -- (11,7) -- (9,4) -- cycle;
\fill[line width=0.7pt,color=qqwwzz,fill=qqwwzz,fill opacity=0.1] (16,4) -- (18,7) -- (20,4) -- cycle;
\draw [line width=0.7pt,color=qqwwzz] (4,5)-- (8,7);
\draw [line width=0.7pt,color=qqwwzz] (8,7)-- (6,3);
\draw [line width=0.7pt,color=qqwwzz] (6,3)-- (4,5);
\draw [line width=0.7pt,color=qqwwzz] (6,3)-- (9,4);
\draw [line width=0.7pt,color=qqwwzz] (9,4)-- (8,7);
\draw [line width=0.7pt,color=qqwwzz] (8,7)-- (6,3);
\draw [line width=0.7pt,color=qqwwzz] (9,4)-- (11,2);
\draw [line width=0.7pt,color=qqwwzz] (11,2)-- (6,3);
\draw [line width=0.7pt,color=qqwwzz] (6,3)-- (9,4);
\draw [line width=0.7pt,color=qqwwzz] (8,7)-- (11,7);
\draw [line width=0.7pt,color=qqwwzz] (11,7)-- (9,4);
\draw [line width=0.7pt,color=qqwwzz] (9,4)-- (8,7);
\draw [line width=0.7pt,color=qqwwzz] (13,4)-- (9,4);
\draw [line width=0.7pt,color=qqwwzz] (9,4)-- (11,2);
\draw [line width=0.7pt,color=qqwwzz] (11,2)-- (13,4);
\draw [line width=0.7pt,color=qqwwzz] (13,4)-- (11,7);
\draw [line width=0.7pt,color=qqwwzz] (11,7)-- (9,4);
\draw [line width=0.7pt,color=qqwwzz] (9,4)-- (13,4);
\draw [line width=0.7pt,color=qqwwzz] (16,4)-- (18,7);
\draw [line width=0.7pt,color=qqwwzz] (18,7)-- (20,4);
\draw [line width=0.7pt,color=qqwwzz] (20,4)-- (16,4);
\draw [line width=0.7pt,dotted] (13,4)-- (16,4);
\draw [line width=0.7pt,dotted] (11,7)-- (18,7);
\draw [line width=0.7pt,color=xfqqff] (16,4)-- (18,5);
\draw [line width=0.7pt,color=xfqqff] (18,5)-- (18,7);
\draw [line width=0.7pt,color=xfqqff] (18,5)-- (20,4);
\draw [-to,line width=0.7pt] (2,1) -- (3,1);
\draw [-to,line width=0.7pt] (2,1) -- (2,2);
\draw [color=qqwwzz](11,5) node[] {$T$};
\draw [color=xfqqff](17.5,5.334) node[] {$T_1$};
\draw [color=xfqqff](18.5,5.334) node[] {$T_2$};
\draw [color=xfqqff](18,4.5) node[] {$T_3$};
\draw (3,1) node[anchor=west] {$x$};
\draw (2,2) node[anchor=south] {$y$};
\begin{scriptsize}
\draw [fill=qqwwzz] (4,5) circle (2.5pt);
\draw [fill=qqwwzz] (6,3) circle (2.5pt);
\draw [fill=qqwwzz] (8,7) circle (2.5pt);
\draw [fill=qqwwzz] (11,2) circle (2.5pt);
\draw [fill=qqwwzz] (13,4) circle (2.5pt);
\draw [fill=qqwwzz] (9,4) circle (2.5pt);
\draw [fill=qqwwzz] (11,7) circle (2.5pt);
\draw [fill=qqwwzz] (16,4) circle (2.5pt);
\draw [fill=qqwwzz] (20,4) circle (2.5pt);
\draw [fill=qqwwzz] (18,7) circle (2.5pt);
\draw [fill=xfqqff] (18,5) circle (2.5pt);
\draw [fill=black] (2,1) circle (2.5pt);
\end{scriptsize}
\end{tikzpicture}
    \caption{Clough--Tocher refinement of an element $T$ into three sub-elements $T_1$, $T_2$ and $T_3$. The sub-elements are not considered as part of a global refinement, but rather define the local space of the macro element $T$.}
    \label{fig:ct}
\end{figure}

In this section we present a finite element that is both, strongly deviatoric, and minimally $\HC{,\surf}$-regular, via a macro element approach. The macro element is derived by the Clough--Tocher refinement of a triangle \cite{Philippe}. Given a triangle $T$, one introduces three new edges between its vertices and the barycentre to split the triangle into three sub-triangles $T = T_1 \cup T_2 \cup T_3$. The refinement is depicted in \cref{fig:ct}, and is not to be understood in the global sense. Rather $T$ is a single finite element.
Now, on each sub-triangle $T_i$ we define the linear polynomial space $\Po^1(\T_i) \otimes \sl(2)$ with the dimension $\dim[\Po^1(\T_i) \otimes \sl(2)] = 3 \cdot 3 = 9$. The complete polynomial space on the macro element is given by $\widetilde{\Y}^1(T) = [\Po^1(\T_1) \otimes \sl(2)] \cup [\Po^1(\T_2) \otimes \sl(2)] \cup [\Po^1(\T_3) \otimes \sl(2)]$ and has the dimension $\dim \widetilde{\Y}^1(T) = 3 \cdot 9 = 27$. Let $\Xi_{ij}$ denote an internal interface (edge) in the macro element between $T_i$ and $T_j$ and $\vb{t}$ its tangent vector, we construct the reduced space
\begin{align}
\label{eq:YT}
    \Y^1(T) = \{ \bm{Y} \in \widetilde{\Y}^1(T) \; | \; \jump{\bm{Y}\vb{t}}|_{\Xi_{ij}} = 0 \quad \forall \, \Xi_{ij} = T_i \cap T_j \} \, ,
\end{align}
such that the space is internally $\HC{,T}$-conforming. Since the polynomial space is sub-element-wise linear, exactly four conditions on each internal edge are required to uphold the conformity, such that $\dim \Y^1(T) = \dim \widetilde{\Y}^1(T) - 3 \cdot 4 = 15$. The internal conditions read 
\begin{align}
    &\int_{\Xi_{ij}}\con{\vb{v}_k \otimes \vb{t}_{ij}}{\bm{Y}} \, \dd \curv = 0 \, , && \vb{v}_k \in [\Po^1(\Xi_{ij})]^2 \, .
\end{align}
With the polynomial space defined, we can construct the macro element.
\begin{definition}[An $\HCd{,\surf}$-conforming macro element] \label{def:macro}
    We define the element via the Ciarlet triplet $\{T,\Y^1(T),S\}$, where the degrees of freedom in $S$ are defined per polytope.
    \begin{itemize}
        \item on each outer edge $e_{ij}$ with tangent vector $\vb{t}_{ij}$ we define
        \begin{align}
            &\int_{\curv_{ij}} \con{\vb{v}_k \otimes \vb{t}_{ij}}{\bm{Y}} \, \dd s \, , && \vb{v}_k \in [\Po^1(\curv_{ij})]^2 \, ,
        \end{align}
        implying four degrees of freedom per edge.
        \item in the cell of $T$ we define the degrees of freedom
        \begin{align}
            &\int_T \con{\bm{V}_i}{\bm{Y}} \, \dd \surf \, , && \bm{V}_i \in \sl(2) \, ,
        \end{align}
        implying three cell functionals.
    \end{itemize}
\end{definition}
\begin{theorem}[Unisolvence of the macro element]
    The construction in \cref{def:macro} defines a unisolvent element for $\HCd{,\surf}$.
\end{theorem}
\begin{proof}
    The proof follows by dimensionality. Observe that $\dim \Y^1(T) = 15$. There are $3 \cdot 4$ edge functionals and $3$ cell functionals, yielding exactly $15$ functionals in total. Therefore, it is sufficient to show that if all degrees of freedom vanish, then the field vanishes as well. Observe that 
    \begin{align}
        \int_T \con{\vb{v}_k}{\Curl \bm{Y}} \, \dd \surf = 0 \quad \forall \, \vb{v}_k \in [\Po^1(T)]^2
    \end{align}
    implies $\Curl \bm{Y} = 0$ in $T$. This becomes apparent by applying the split over the sub-triangles, which is possible due to interelement $\HC{}$-conformity \eqref{eq:YT},
    \begin{align}
        \int_T \con{\vb{v}_k}{\Curl \bm{Y}} \, \dd \surf = \sum_{i=1}^3 \int_{T_i} \con{\vb{v}_k}{\Curl \bm{Y}} \, \dd \surf = \sum_{i=1}^3 |T_i| \con{\vb{v}_k}{\Curl \bm{Y}}\at_{\vb{x}_c} \, ,  
    \end{align}
    and applying the Gaussian one-point quadrature rule, which is exact for linear polynomials. 
    Now, by Green's formula we find
    \begin{align}
        \int_T \con{\vb{v}_k}{\Curl \bm{Y}} \, \dd \surf = \int_{\partial T} \con{\vb{v}_k}{ \bm{Y} \vb{t} } \, \dd \curv + \int_T \con{\dev \D^\perp \vb{v}_k}{ \bm{Y}} \, \dd \surf \, ,
    \end{align}
    such that if all the degrees of freedom are zero, then the equation yields zero and we have $\Curl \bm{Y}|_T = 0$. By applying the Stokes theorem sub-element-wise we find
    \begin{align}
        \int_{T_i} \Curl \bm{Y} \, \dd \surf = \underbrace{\int_{\curv_j} \bm{Y} \vb{t} \, \dd \curv}_{c_j} + \underbrace{\int_{\curv_k} \bm{Y} \vb{t} \, \dd \curv}_{c_k} = 0 \, , 
    \end{align}
    where $\curv_j$ and $\curv_k$ represent the internal interface-edges of each sub-triangle. Repeating the procedure for each sub-element while denoting the integral terms on each sub-element $c_q$, we find the system of equations
    \begin{align}
        &c_1 + c_2 = 0 \, , && c_2 + c_3 = 0 \, , && c_3 + c_1 = 0 \, ,
    \end{align}
    whose solution is $c_1 = c_2 = c_3 = 0$. Consequently, we have by Green's formula on each sub-triangle 
    \begin{align}
        \underbrace{\int_{T_i} \con{\vb{v}_k}{\Curl \bm{Y}} \, \dd \surf}_{= 0} = \underbrace{\int_{\partial T_{i}} \con{\vb{v}_k}{ \bm{Y} \vb{t} } \, \dd \curv}_{= 0} + \int_{T_i} \con{\dev\D^\perp \vb{v}_k}{ \bm{Y}} \, \dd \surf = 0 \quad \forall \, \vb{v}_k \in [\Po^1(T_i)]^2 \, ,
    \end{align}
    implying $\bm{Y}|_T = 0$.
\end{proof}

\subsection{A weakly deviatoric formulation} \label{sec:weak}
The construction in \cref{sec:dev} is strongly deviatoric but not minimally $\HC{,\surf}$-regular due to the higher continuity at vertices.  
The lack of minimal regularity may lead to an incorrect imposition of boundary conditions. Further, in the case of jumping material coefficients between elements, the lack of minimal regularity can lead to errors in the approximation. 
While the construction in \cref{sec:macro} is strongly deviatoric and minimally regular, by its nature as a macro element it does not lend itself to a straight-forward implementation, such that we do not consider it further in this current work.
Alternatively, we restore minimal regularity 
by abandoning strong tracelessness and then re-imposing it weakly with a mixed formulation
\begin{align}
    a(\{\delta\ud,\delta\devP,\delta\sphP\},\{\ud,\devP,\sphP\}) + \con{\delta\devP}{\bm{V}}_{\Le} &= l(\{\delta\ud,\delta\devP,\delta\sphP\}) \, , \notag \\
    \con{\delta\bm{V}}{\devP}_{\Le} &= 0 \, ,
\end{align}
where $\delta \bm{V},\bm{V} \in [\Le(\surf) \cdot \one]$ are discontinuous spherical tensors and the deviatoric microdistortion is now $\devP \in \HC{,\surf}$. Clearly, the constraint given by $\con{\delta \bm{V}}{\devP}_{\Le} = 0$ implies tracelessness of the microdistortion field $\tr \devP = 0$. The discretisation and stability of the formulation follows analogously to the weak elasticity complex \cite{arnold_mixed_2007}, such that $\bm{V} \in [\DG^{p-1}(\surf) \cdot \one]$, and $\devP \in [\R^2 \otimes \Nedtwo^{p}(\surf)]$ is a row-wise N\'ed\'elec element of the second type \cite{sky_polytopal_2022,sky_higher_2023,Ned2,haubold2023high}.  

\subsection{Convergence estimates}
The primal formulation \cref{eq:primal_form} is shown to be well-posed by the Lax--Milgram theorem, such that we can readily apply Cea's lemma to find an a priori error estimate when using the strongly deviatoric elements $\Y^{p}(\surf)$.
\begin{theorem}[Convergence of the primal form] \label{thm:conv}
    Let the exact solution be $\{\ud,\devP,\sphP\} \in [\H^{p+1}(\surf)]^2 \times [\H^{p+1}(\surf)\otimes\sl(2)] \times [\H^{p}(\surf) \cdot \one]$, then the discrete solution $\{\ud^h,\devP^h,\sphP^h\} \in \X_h(\surf) = [\CG^{p+1}(\surf)]^2 \times \Y^{p}(\surf) \times [\DG^{p}(\surf) \cdot \one] \subset \X(\surf)$ satisfies the convergence rate 
    \begin{align}
    \norm{\{\ud,\devP,\sphP\} - \{\ud^h,\devP^h,\sphP^h\}}_\X & \leq c\,h^{p}(|\ud|_{\H^{p+1}} + |\devP|_{\H^{p+1}} + |\Curl\devP|_{\H^{p}} + |\sphP|_{\H^p}) \, ,
    \end{align}
    under a uniform triangulation. 
\end{theorem}
\begin{proof}
    Let $\Pi_g^{p+1}:\Hone(\surf) \to \CG^{p+1}(\surf)$, $\Pi_o^{p}:[\Le(\surf)\cdot \one] \to [\DG^{p}(\surf)\cdot \one]$ and $\Pi_y^{p}:\HCd{,\surf}\to\Y^{p}(\surf)$, there holds
    \begin{align}
    \norm{\{\ud,\devP,\sphP\} - \{\ud^h,\devP^h,\sphP^h\}}_\X^2 & \leq c \inf_{\{\delta\ud^h,\delta\devP^h,\delta\sphP^h\} \in \X_h} \norm{\{\ud,\devP,\sphP\} - \{\delta\ud^h,\delta\devP^h,\delta\sphP^h\}}_\X^2 \notag \\
    & \leq c \, (\norm{\ud - \Pi_g^{p+1}\ud}_{\Hone}^2 + \norm{\devP - \Pi_y^{p}\devP}_{\HC{}}^2 + \norm{\sphP - \Pi_o^{p}\sphP}_{\Le}^2) \notag \\
    & \leq c \, (h^{2p}|\ud|_{\H^{p+1}}^2 + \norm{\devP - \Pi_y^{p}\devP}_{\Le}^2 + \norm{\Curl\devP - \Curl\Pi_y^{p}\devP}_{\Le}^2 + h^{2p}|\sphP|_{\H^{p}}^2) \notag \\
    & \leq c \, (h^{2p}|\ud|_{\H^{p+1}}^2 + h^{2(p+1)}|\devP|_{\H^{p+1}}^2 + h^{2p}|\Curl\devP|_{\H^{p}}^2 + h^{2p}|\sphP|_{\H^{p}}^2) \notag \\
    & \leq c \, h^{2p} (|\ud|_{\H^{p+1}}^2 + |\devP|_{\H^{p+1}}^2 + |\Curl\devP|_{\H^{p}}^2 + |\sphP|_{\H^{p}}^2) \, ,\label{eq:conv_proof1}
    \end{align}
    where we used Cea's lemma, standard continuous and discontinuous Lagrange interpolants, and finally the estimates from \cref{eq:inty}. 
\end{proof}

Since the strongly deviatoric element satisfies a commuting grad-curl sub-complex, we can discretise the mixed formulation as $\{\ud,\devP,\sphP, \hP\} \in \X_h(\surf) \times [\CG^{p-1}(\surf)]^2$, with $\X_h$ as in \cref{thm:conv} and $p \geq 3$. 
\begin{theorem}[Convergence of the mixed form] \label{thm:conv2}
    Let the exact solution be $\{\ud,\devP,\sphP,\hP\} \in [\H^{p+1}(\surf)]^2 \times [\H^{p+1}(\surf)\otimes\sl(2)] \times [\H^{p}(\surf) \cdot \one] \times [\H^{p}(\surf)]^2$, then the discrete solution $\{\ud^h,\devP^h,\sphP^h,\hP^h\} \in \X_h(\surf) \times \Z_h(\surf) = \X_h(\surf) \times [\DG^{p-1}(\surf)]^2  \subset \X(\surf) \times \Z(\surf) = \X(\surf) \times [\Le(\surf)]^2$ satisfies the convergence rate 
    \begin{align}
    \norm{\{\ud,\devP,\sphP,\hP\} - \{\ud^h,\devP^h,\sphP^h,\hP^h\}}_{\X\times\Z} & \leq c \, h^{p}(|\ud|_{\H^{p+1}} + |\devP|_{\H^{p+1}} + |\Curl\devP|_{\H^{p}} + |\sphP|_{\H^{p}} + |\hP|_{\H^{p}}) \, ,
    \end{align}
    under a uniform triangulation, where the constant $c$ does not depend on $\Lc$. 
\end{theorem}
\begin{proof}
    Let $\W(\surf) = \X(\surf) \times \Z(\surf)$ and $\W_h(\surf) = \X_h(\surf) \times \Z_h(\surf)$, the proof follows by extending the previous convergence estimate \cref{eq:conv_proof1}
    \begin{align}
    \norm{\{\ud,\devP,\sphP,\hP\} - \{\ud^h,\devP^h,\sphP^h,\hP^h\}}_{\W}^2 & \leq c \inf_{\{\delta\ud^h,\delta\devP^h,\delta\sphP^h,\delta\hP^h\} \in {\W_h}} \norm{\{\ud,\devP,\sphP,\hP\} - \{\delta\ud^h,\delta\devP^h,\delta\sphP^h,\delta\hP^h\}}_{\W}^2 
    \notag \\
    & \leq c  \,  (\norm{\ud - \Pi_g^{p+1}\ud}_{\Hone}^2 + \norm{\devP - \Pi_y^{p}\devP}_{\HC{}}^2 + \norm{\sphP - \Pi_o^{p}\sphP}_{\Le}^2 + \norm{\hP - \Pi_o^{p-1}\hP}_{\Le}^2)
    \notag \\
    & \leq c \,  (h^{2p}|\ud|_{\H^{p+1}}^2 + h^{2(p+1)}|\devP|_{\H^{p+1}}^2 + h^{2p}|\Curl\devP|_{\H^{p}}^2 + h^{2p}|\sphP|_{\H^{p}}^2 + h^{2p}|\hP|_{\H^{p}}^2) 
    \notag \\
    & \leq c \,  h^{2p} (|\ud|_{\H^{p+1}}^2 + |\devP|_{\H^{p+1}}^2 + |\Curl\devP|_{\H^{p}}^2 + |\sphP|_{\H^{p}}^2 + |\hP|_{\H^{p}}^2) \, ,
    \end{align}
    to account for the hyper-stress field $\hP$.
\end{proof}
\begin{remark}[Convergence of the weakly deviatoric formulation]
    Analogously to \cref{thm:conv2}, we find the same convergence estimate for the primal weakly deviatoric formulation in \cref{sec:weak}. This becomes apparent when replacing $\Pi_y^p \devP$ with the row-wise application of the second-type N\'ed\'elec interpolant $\Pi_c^p \devP$, and $\Pi_o^{p-1}\hP$ with $\Pi_o^{p-1} \bm{V}$. This is because the second-type N\'ed\'elec interpolant yields the same convergence estimates.    
\end{remark}

\section{Numerical benchmarks}
In the following we consider two examples of the newly introduced planar relaxed micromorphic model. The first example demonstrates the capacity of the formulation and our finite elements to correctly capture discontinuous dilatation fields, and compares it with the full Curl planar relaxed micromorphic model. In the second example we verify optimal convergence and stability of our finite element discretisations, as well as demonstrate the homogeneous behaviour of the model for $\Lc \to 0$. 

\subsection{A discontinuous planar dilatation field}
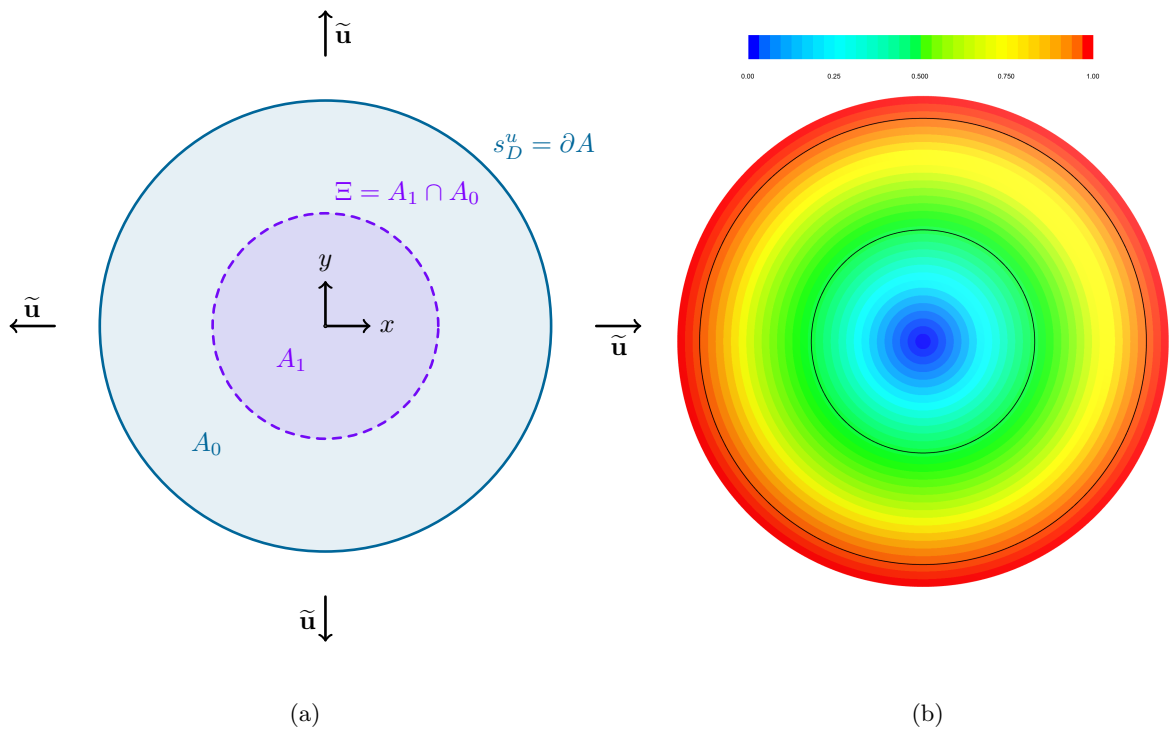
\begin{figure}
		\centering
		\begin{subfigure}{0.48\linewidth}
			\centering
		\definecolor{qqwwzz}{rgb}{0,0.4,0.6}
\definecolor{xfqqff}{rgb}{0.4980392156862745,0,1}
\definecolor{uuuuuu}{rgb}{0.26666666666666666,0.26666666666666666,0.26666666666666666}
\begin{tikzpicture}[scale = 0.3, line cap=round,line join=round,>=triangle 45,x=1cm,y=1cm]
\clip(-14.5,-14.5) rectangle (14.5,14.5);
\draw [line width=1pt,dashed,color=xfqqff,fill=xfqqff,fill opacity=0.1] (0,0) circle (5cm);
\draw [line width=1pt,color=qqwwzz,fill=qqwwzz,fill opacity=0.1] (0,0) circle (10cm);
\draw [-to,line width=1pt] (0,0) -- (2,0);
\draw [-to,line width=1pt] (0,0) -- (0,2);
\draw [color=black](2,0) node[anchor=west] {$x$};
\draw [color=black](0,2) node[anchor=south] {$y$};
\draw [color=qqwwzz](7,7) node[anchor=south west] {$s_D^u = \partial A$};
\draw [color=qqwwzz](-5.25,-5.25) node[] {$A_0$};
\draw [color=xfqqff](-1.5,-1.5) node[] {$A_1$};
\draw [color=xfqqff](0,5) node[anchor=south west] {$\Xi = A_1 \cap A_0$};
\draw [-to,line width=1pt] (0,-12) -- (0,-14);
\draw [-to,line width=1pt] (12,0) -- (14,0);
\draw [-to,line width=1pt] (0,12) -- (0,14);
\draw [-to,line width=1pt] (-12,0) -- (-14,0);
\draw (13,0) node[anchor=north] {$\widetilde{\mathbf{u}}$};
\draw (0,13) node[anchor=west] {$\widetilde{\mathbf{u}}$};
\draw (-13,0) node[anchor=south] {$\widetilde{\mathbf{u}}$};
\draw (0,-13) node[anchor=east] {$\widetilde{\mathbf{u}}$};
\begin{scriptsize}
\draw [fill=uuuuuu] (0,0) circle (2pt);
\end{scriptsize}
\end{tikzpicture}
			\caption{}
		\end{subfigure}
	    \begin{subfigure}{0.48\linewidth}
	    	\centering
	    \includegraphics[width=0.8\linewidth]{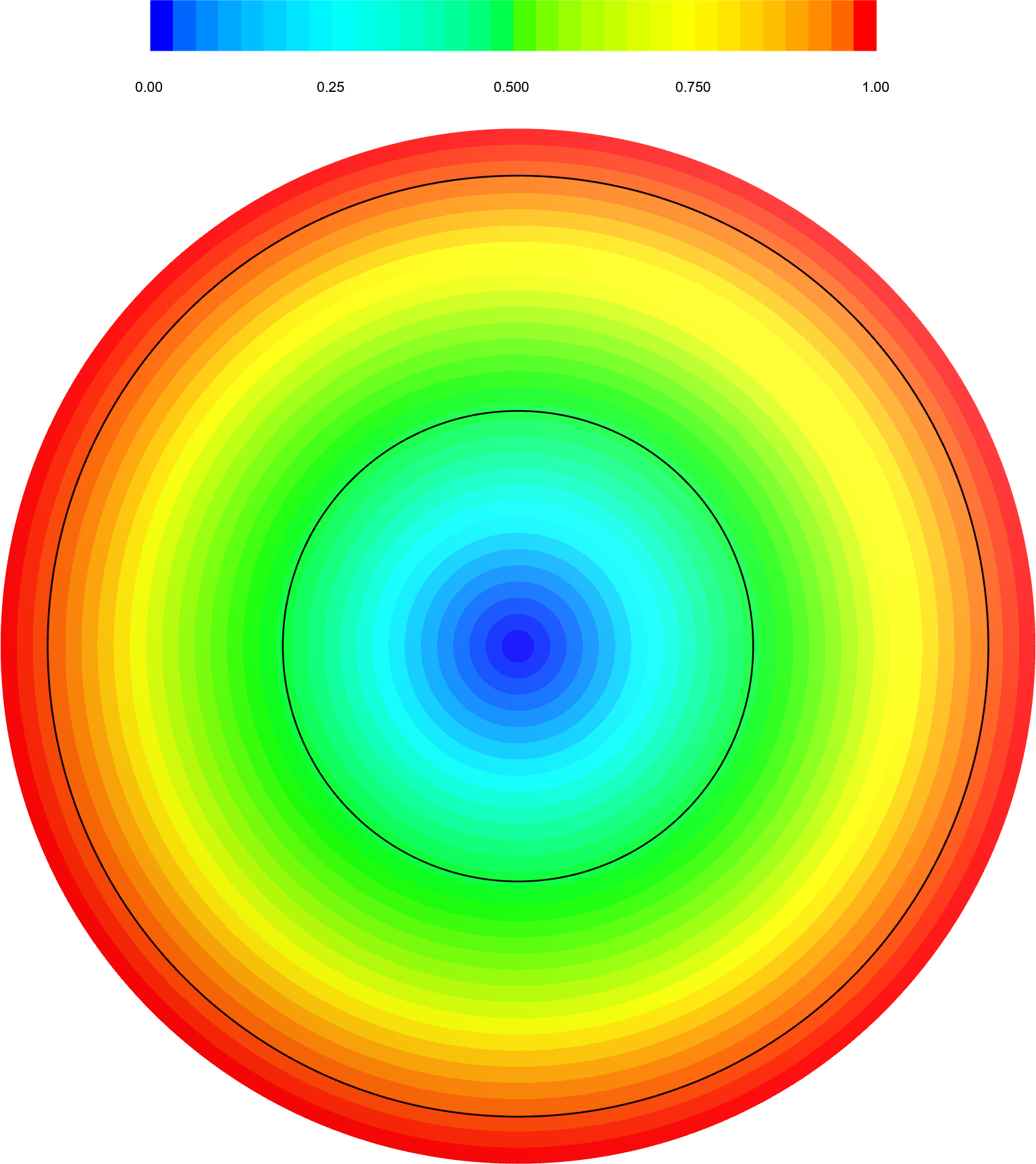}
     \vspace{1.3cm}
	    	\caption{}
	    \end{subfigure}
		\caption{Domain $\surf = \surf_0 \cup \surf_1$ composed of two materials that jump at the interface (a). The interface between the inner $\surf_1$ and outer $\surf_0$ domains is given by $\Xi$. On the outer boundary $\partial \surf_0 = \curv_D^u$ we impose the Dirichlet boundary condition for the displacement field $\vb{u}$, which is depicted in (b) for a mesh of $2636$ elements.}
\label{fig:dom}
\end{figure}

In the following we provide an example for the occurrence of a jumping dilatation field under uniform expansion when the Lam\'e moduli differ between interfacing materials.
We define the global material constants
\begin{align}
    &\lamma = 115.4 \, , && \muma = 76.9 \, , && \muc = 0 \, , && \Lc = 1 \, ,
\end{align}
and a circular domain composed of two materials $\overline{\surf} = \overline{\surf}_{0} \cup \overline{\surf}_{1}$, such that $\surf_0 = \{(x,y) \in (0,10)^2 \; | \; x^2 + y^2 < 100\}$ and $\surf_1 = \{(x,y) \in (0,5)^2 \; | \; x^2 + y^2 < 25\}$, see \cref{fig:dom}.
On the outer domain $\surf_0$ we set the micro parameters to
\begin{align}
    & \lammi^0 = 10 \, \lamma \, , && \mumi^0 = 10 \, \muma \, , 
\end{align}
by scaling the macro parameters. On the inner domain $\surf_1$ we define
\begin{align}
    & \lammi^1 = 1000 \, \lamma \, , && \mumi^1 = 1000 \, \muma \, , 
\end{align}
implying that the inner domain is composed of much stiffer material than the outer domain.
The remaining material constants for each material are retrieved via the planar homogenisation formula \cite{GOURGIOTIS2024112700}
\begin{align}
    \mue = \dfrac{\mumi\, \muma}{\mumi- \muma} \, , && \mue +  \lambda_\mathrm{e} = \dfrac{(\mumi + \lambda_\mathrm{micro})(\muma + \lambda_\mathrm{macro})}{(\mumi + \lambda_\mathrm{micro})-(\muma + \lambda_\mathrm{macro})} \, .
    \label{eq:homogen}
\end{align}
On the boundary of the domain $\partial \surf$ we impose the Dirichlet boundary condition 
\begin{align}
    &\vb{u} \at_{\curv_D} = \widetilde{\vb{u}} = \dfrac{1}{10} \begin{bmatrix}
        x \\ y
    \end{bmatrix} \, , && \curv_D^u = \partial \surf \, .
\end{align}
The boundary condition of the deviatoric microdistortion field $\devP$ is a complete Neumann boundary with no couple tractions. In the entire domain, the forces $\vb{f}$ and couple forces $\bm{M}$ are set to zero.

\begin{figure}
		\centering
		\begin{subfigure}{0.48\linewidth}
			\centering
		\begin{tikzpicture}[scale = 0.9]
    			\definecolor{asl}{rgb}{0.4980392156862745,0.,1.}
    			\definecolor{asb}{rgb}{0.,0.4,0.6}
    			\begin{semilogxaxis}[
    				/pgf/number format/1000 sep={},
    				axis lines = left,
    				xlabel={Number of elements},
    				ylabel={Lebesgue norm},
    				xmin=50, xmax=3e3,
    				ymin=0.01, ymax=0.25,
    				xtick={1e1,1e2,1e3,1e4,1e5,1e6},
    				ytick={0.02, 0.07, 0.13, 0.18, 0.23},
    				legend style={at={(0.95,0.5)},anchor= east},
    				ymajorgrids=true,
    				grid style=dotted,
    				]
    				\addplot[color=asl, mark=pentagon] coordinates {
    					( 54 , 0.2170839902975513 )
                            ( 168 , 0.21708399193998515 )
                            ( 602 , 0.21708399430577785 )
                            ( 2636 , 0.21708399434301937 )
    				};
    				\addlegendentry{$\norm{\sphP}_{\Le}$}
    				
    				\addplot[color=asb, mark=triangle] coordinates {
                            ( 54 , 0.21517476207817499 )
                            ( 168 , 0.2158304521021457 )
                            ( 602 , 0.21588354722152434 )
                            ( 2636 , 0.21591991655666226 )
    				};
    				\addlegendentry{$\norm{\sph \Pm}_{\Le}$}

                      \addplot[color=blue, mark=diamond] coordinates {
                            ( 54 , 0.028287213187060217 )
                            ( 168 , 0.09669869536747962 )
                            ( 602 , 0.01823314539733175 )
                            ( 2636 , 0.04733260935781584 )
    				};
    				\addlegendentry{$\norm{\dev \Pm}_{\Le}$}
    				
    			\end{semilogxaxis}


\end{tikzpicture}
			\caption{}
		\end{subfigure}
	    \begin{subfigure}{0.48\linewidth}
	    	\centering
	    \includegraphics[width=0.8\linewidth]{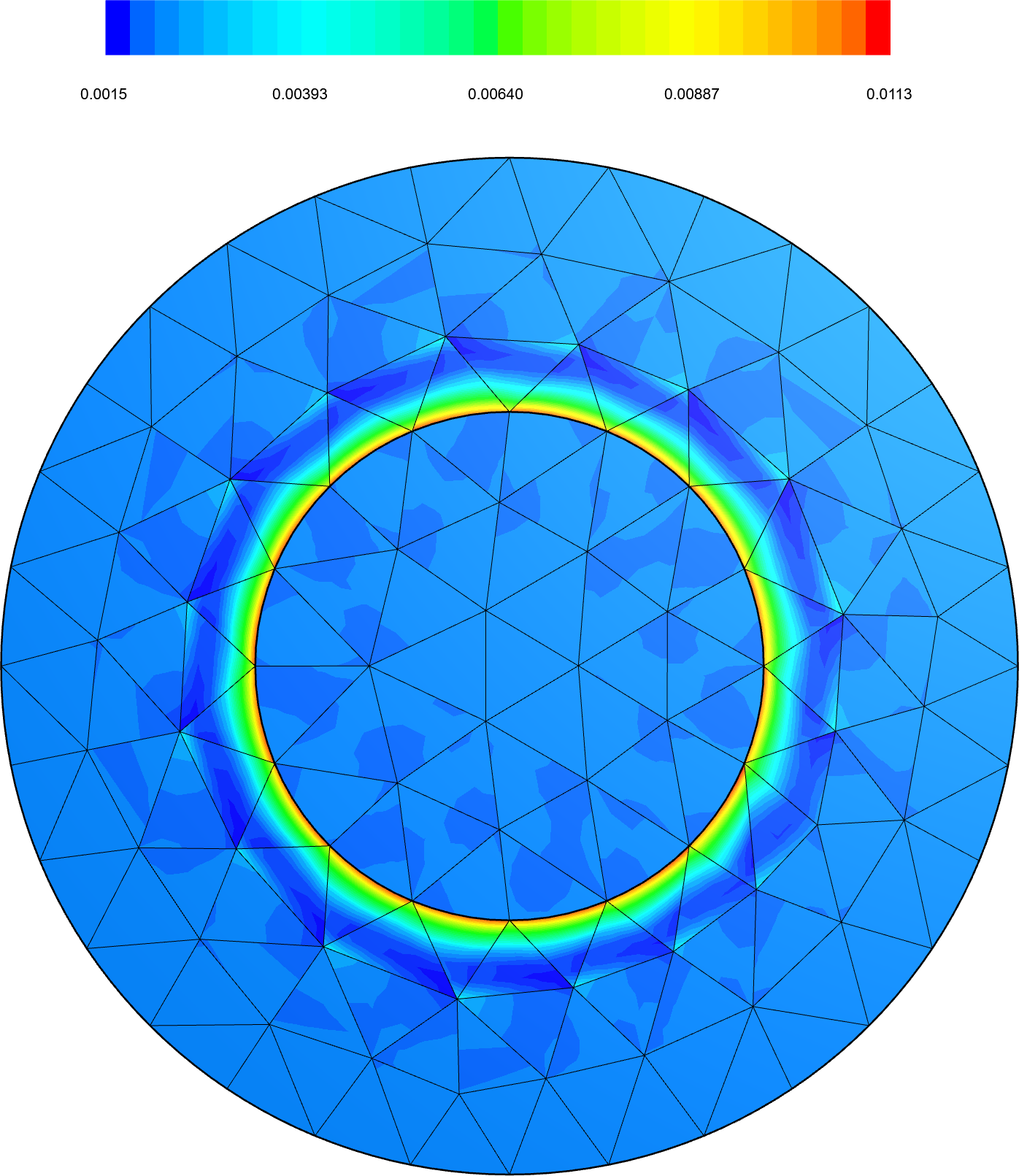}
	    	\caption{}
	    \end{subfigure}
        \begin{subfigure}{0.48\linewidth}
	    	\centering
      \vspace{0.1cm}
	    \includegraphics[width=0.8\linewidth]{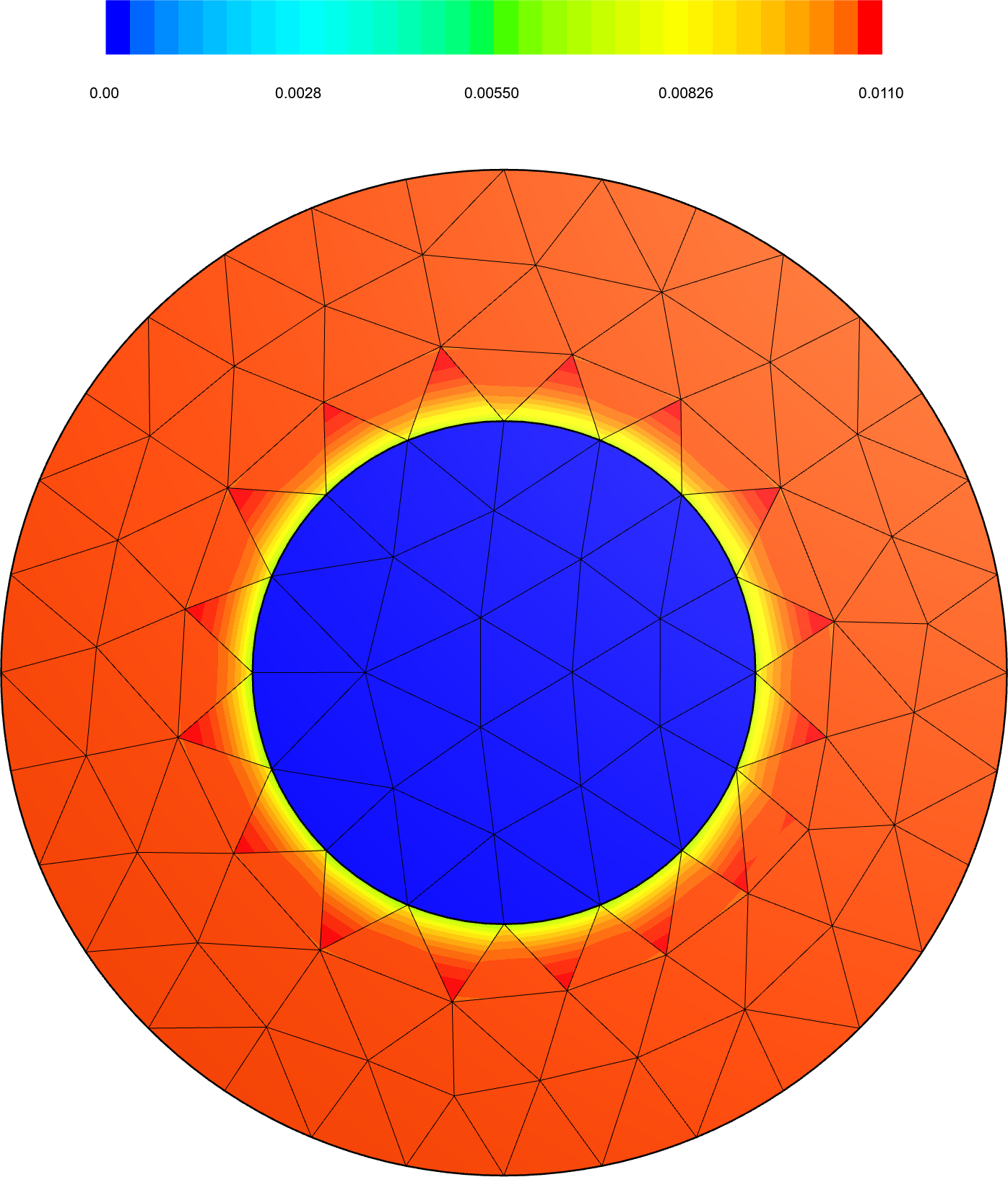}
	    	\caption{}
	    \end{subfigure} 
        \begin{subfigure}{0.48\linewidth}
	    	\centering
      \vspace{0.1cm}
	    \includegraphics[width=0.8\linewidth]{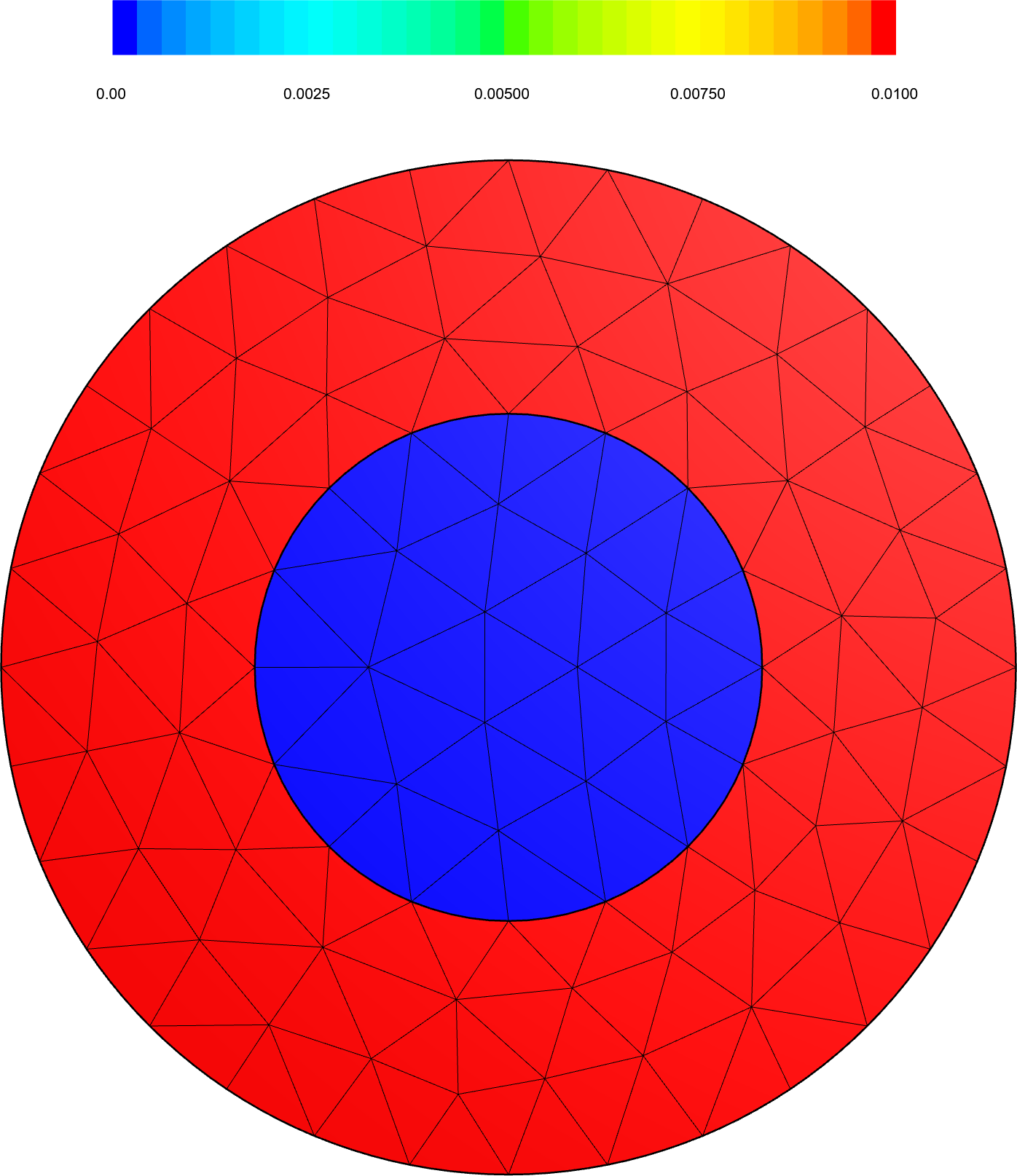}
	    	\caption{}
	    \end{subfigure}
		\caption{Norms of the microdistortion fields in the full Curl formulation and the newly presented formulation under mesh refinement (a). The fields $\norm{\dev \Pm}$ and $\sph \Pm$ are depicted in (b) and (c) for the full Curl formulation, respectively. The field $\devP$ is zero across the domain in new formulation, while the field $\sphP$ is depicted in (d).}
		\label{fig:conv}
\end{figure}
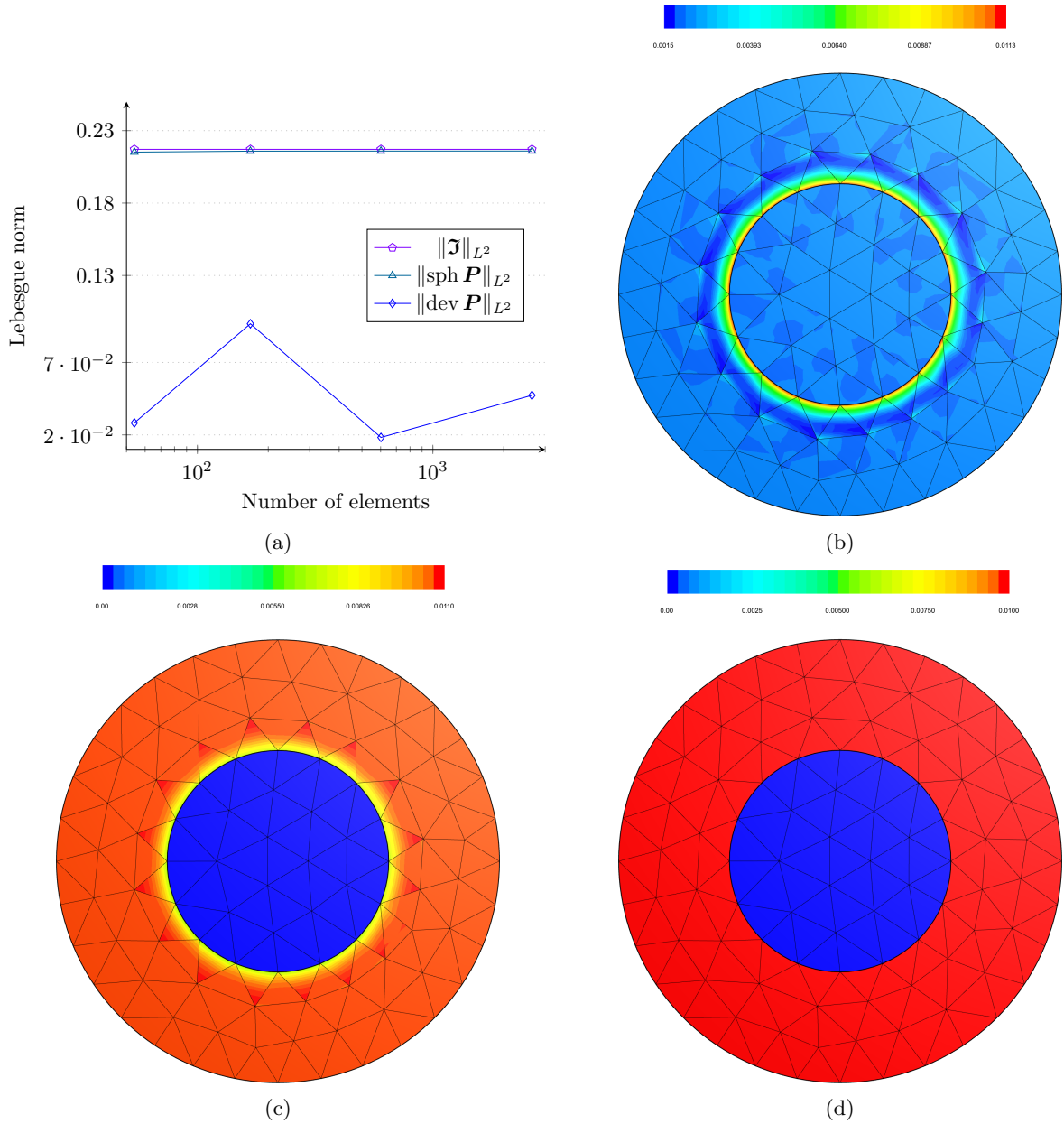

We compare the new planar formulation in \cref{eq:redrmm} with the full Curl formulation of the relaxed micromorphic model   
\begin{align}
    &\int_\surf \con{\Ce\sym(\D \delta\ud - \delta\Pm)}{\sym(\D \ud - \Pm)} + \con{\Cc\skw(\D \delta\ud - \delta\Pm)}{\skw(\D \ud - \Pm)} + \con{\Cmic\sym\delta\Pm}{\sym\Pm} \notag \\ &\quad + \muma \Lc^2 \con{\Curl \delta \Pm}{\Curl \Pm} \, \dd \surf = \int_\body \con{\delta\ud}{\vb{f}} + \con{\delta\Pm}{\bm{M}} \, \dd \surf \, .
\end{align}
In both cases we use quartic polynomial approximations $p = 4$.
For the new planar formulation we use the deviatoric element from \cref{sec:dev} denoted as $\mathcal{Y}^p(\surf)$ in the discretisation $\{\ud, \devP, \sphP\} \in [\CG^4(\surf)]^2 \times \mathcal{Y}^3(\surf) \times [\DG^3(\surf) \cdot \one] \subset [\Hone(\surf)]^2 \times \HCd{,\surf} \times [\Le(\surf) \cdot \one]$. For the full Curl formulation we use $\{\ud, \Pm\} \in [\CG^4(\surf)]^2 \times [\R^2 \otimes\Nedtwo^3(\surf) ] \subset [\Hone(\surf)]^2 \times \HC{,\surf}$, where $[\R^2 \otimes\Nedtwo^3(\surf)]$ are row-wise cubic N\'ed\'elec elements of the second type \cite{sky_polytopal_2022,sky_higher_2023,SKY2022115298,SkyPamm}. The change in the norms of the isochoric decomposition of the microdistortion field of the respective formulations for meshes with $54$, $168$, $602$ and $2636$ quartic elements is measured in \cref{fig:conv} (a).
Clearly, mesh refinement does not fix the underlying incapacity of the full Curl model to capture a jumping dilatation field. Consequently, the norm of the deviatoric part of the microdistortion $\norm{\dev \Pm}$ never vanishes for the full Curl formulation regardless of refinement. In contrast, the modified form can perfectly accommodate a jumping dilatation field, such that $\norm{\devP}$ is numerically zero. In \cref{fig:conv} (b) we depict the norm of non-vanishing deviatoric part of the microdistortion $\norm{\dev \Pm}$ of the full Curl formulation. Note that since the field $\devP$ completely vanishes we do not depict it. In \cref{fig:conv} (c) we observe that the full Curl formulation attempts to compensate for the jump in the material coefficients with a sharp transition in $\sph \Pm$. In contrast, the new formulation allows to exactly capture the jump. The latter is shown in \cref{fig:conv} (d), where we depict the field $\sphP$.

\begin{figure}
		\centering
		\begin{subfigure}{0.48\linewidth}
			\centering
		\includegraphics[width=0.8\linewidth]{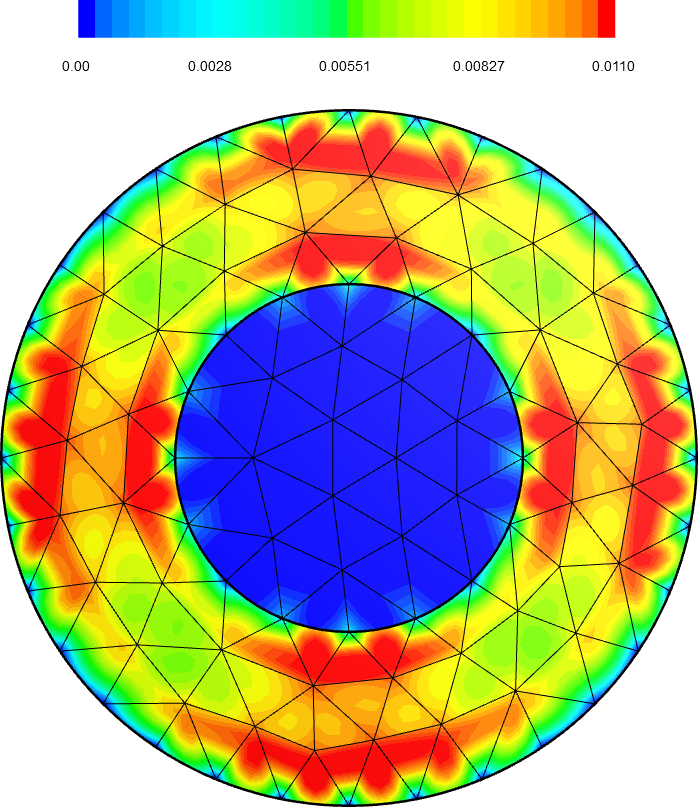}
			\caption{}
		\end{subfigure}
	    \begin{subfigure}{0.48\linewidth}
	    	\centering
	    \includegraphics[width=0.8\linewidth]{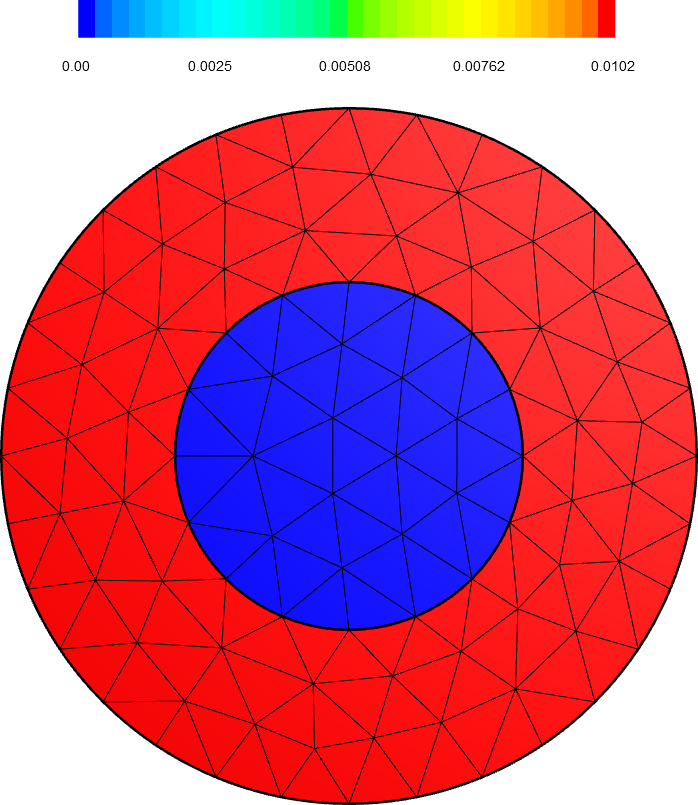}
	    	\caption{}
	    \end{subfigure}
        \begin{subfigure}{0.48\linewidth}
	    	\centering
      \vspace{0.1cm}
	    \includegraphics[width=0.8\linewidth]{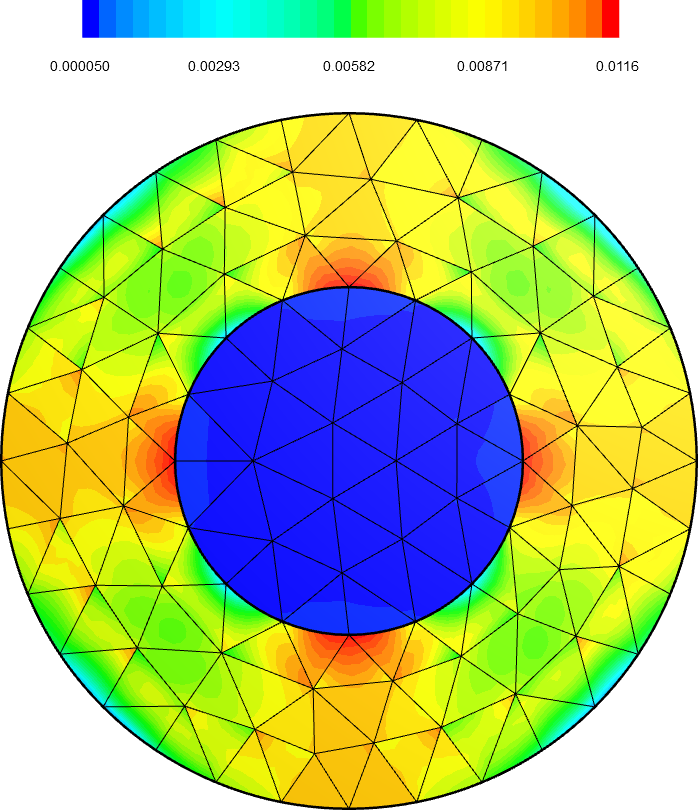}
	    	\caption{}
	    \end{subfigure} 
        \begin{subfigure}{0.48\linewidth}
	    	\centering
      \vspace{0.1cm}
	    \includegraphics[width=0.8\linewidth]{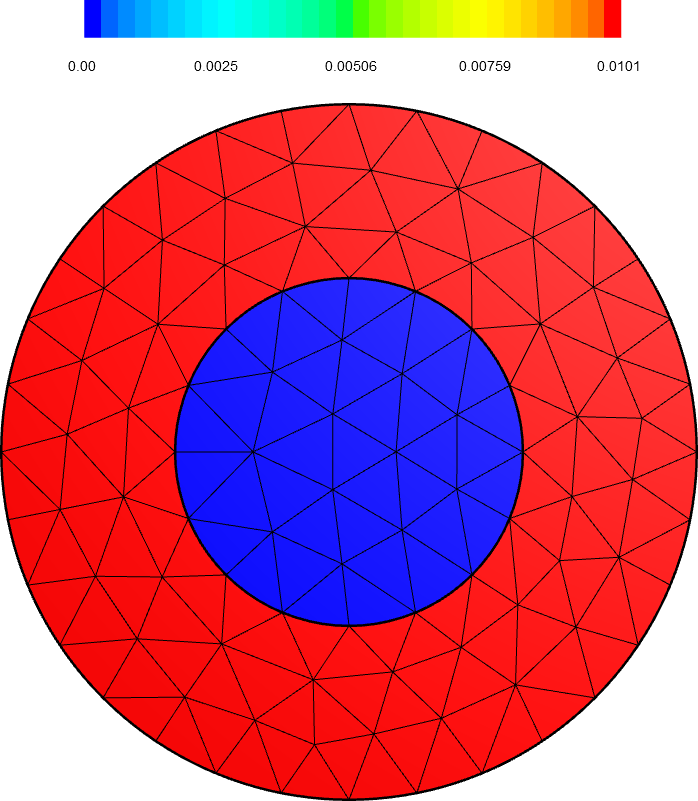}
	    	\caption{}
	    \end{subfigure}
		\caption{Norm of the strongly deviatoric microdistortion field $\norm{\devP}$ (a), and the spherical $\sphP$ (b) microdistortion field. In contrast, in (c) the norm of the deviatoric microdistortion field is defined as $\norm{\dev \devP}$, since the field is only weakly traceless. Correspondingly, in (d) the spherical microdistortion field is $\sphP + \sph \devP$.}
		\label{fig:qcomp}
\end{figure}

Taking the same domain and material parameters, we slightly modify the boundary conditions to
\begin{align}
    &\vb{u} \at_{\curv_D} = \widetilde{\vb{u}} = \dfrac{1}{10} \begin{bmatrix}
        x \\ y + x
    \end{bmatrix} \, ,  && \curv_D = \partial \surf \, ,
\end{align}
as well as impose zero on the part of the deviatoric microdistortion field $\devP$ we control. Consequently, the deformation now includes shear and the deviatoric microdistortion field does not vanish. 
Further, due to the differing shear moduli, the deviatoric field must jump between the inner and outer materials.
Since we compare here the strongly deviatoric formulation (\cref{sec:dev}) with the weakly deviatoric formulation (\cref{sec:weak}), the boundary conditions on $\devP$ are not equivalent, which makes a quantitative comparison difficult. As such, we compare the behaviour of the strongly deviatoric and weakly deviatoric element qualitatively in \cref{fig:qcomp}.  
We observe that both formulations correctly account for the jumping dilatation field. However, only the weakly deviatoric formulation, where the element is minimally $\HC{,\surf}$-regular, is capable of capturing the jump in the shear-distortion.  

\subsection{Convergence and limits of $\Lc$}
\begin{figure}
    	\centering
    	\begin{subfigure}{0.48\linewidth}
    		\centering
    		\begin{tikzpicture}[scale = 0.9]
    			\definecolor{asl}{rgb}{0.4980392156862745,0.,1.}
    			\definecolor{asb}{rgb}{0.,0.4,0.6}
    			\begin{loglogaxis}[
    				/pgf/number format/1000 sep={},
    				axis lines = left,
    				xlabel={Degrees of freedom},
    				ylabel={Error},
    				xmin=900, xmax=40000,
    				ymin=0.75e-6, ymax=0.05,
    				xtick={1e2,1e3,1e4,1e5},
    				ytick={1e-6,1e-4,1e-2},
    				legend style={at={(0.95,0.95)},anchor= north east},
    				ymajorgrids=true,
    				grid style=dotted,
    				]
    				\addplot[color=asb, mark=diamond] coordinates {
                            ( 1485 , 0.004880769884549849 )
                            ( 5685 , 0.0007531063058481982 )
                            ( 12605 , 0.00010660380103566646 )
                            ( 17085 , 5.014451291410993e-05 )
                            ( 22245 , 2.600553650941184e-05 )
    				};
    				\addlegendentry{$\norm{\widetilde{\ud} - \ud^h}_{\Le}$}

                    \addplot[color=blue, mark=triangle] coordinates {
                            ( 1485 , 0.0003949850891879821 )
                            ( 5685 , 7.134875819503008e-05 )
                            ( 12605 , 9.33368385614957e-06 )
                            ( 17085 , 4.2205937659323444e-06 )
                            ( 22245 , 2.1107625834573216e-06 )
    				};
    				\addlegendentry{$\norm{\widetilde{\devP} - \devP^h}_{\Le}$}
    				
    				\addplot[color=asl, mark=square] coordinates {
        ( 1485 , 0.00016902136346478374 )
        ( 5685 , 4.190326648221097e-05 )
        ( 12605 , 7.281008609392433e-06 )
        ( 17085 , 3.75423690747589e-06 )
        ( 22245 , 2.1208602544016252e-06 )
    				};
    				\addlegendentry{$\norm{\widetilde{\sphP} - \sphP^h}_{\Le}$}
    				
    				\addplot[dashed,color=black, mark=none]
    				coordinates {
    					(6000, 9e-5)
    					(20000, 4.436552715791846e-06)
    				};    		
                    \addplot[dashed,color=black, mark=none]
    				coordinates {
    					(6000, 2e-5)
    					(20000, 1.8e-6)
    				};    		
    			\end{loglogaxis}
    			\draw (4.45,1.45) 
    			node[anchor=south west]{$\mathcal{O}(h^{5})$};
                \draw (5.05,1) 
    			node[anchor=north east]{$\mathcal{O}(h^{4})$};
    		\end{tikzpicture}
    		\caption{}
    	\end{subfigure}
    	\begin{subfigure}{0.48\linewidth}
    		\centering
    		\begin{tikzpicture}[scale = 0.9]
    			\definecolor{asl}{rgb}{0.4980392156862745,0.,1.}
    			\definecolor{asb}{rgb}{0.,0.4,0.6}
    			\begin{loglogaxis}[
    				/pgf/number format/1000 sep={},
    				axis lines = left,
    				xlabel={Degrees of freedom},
    				ylabel={Error},
    				xmin=900, xmax=40000,
    				ymin=0.75e-6, ymax=0.05,
    				xtick={1e2,1e3,1e4,1e5},
    				ytick={1e-6,1e-4,1e-2},
    				legend style={at={(0.05,0.05)},anchor= south west},
    				ymajorgrids=true,
    				grid style=dotted,
    				]
    				\addplot[color=asb, mark=diamond] coordinates {
                            ( 1871 , 0.004880769884549851 )
                            ( 7223 , 0.0007531063058481986 )
                            ( 16063 , 0.00010660380103566629 )
                            ( 21791 , 5.014451291410979e-05 )
                            ( 28391 , 2.6005536509411885e-05 )
    				};
    				\addlegendentry{$\norm{\widetilde{\ud} - \ud^h}_{\Le}$}

                    \addplot[color=blue, mark=triangle] coordinates {
                            ( 1871 , 0.00039498508918799953 )
                            ( 7223 , 7.134875819503697e-05 )
                            ( 16063 , 9.33368385615042e-06 )
                            ( 21791 , 4.2205937659327375e-06 )
                            ( 28391 , 2.1107625834574652e-06 )
    				};
    				\addlegendentry{$\norm{\widetilde{\devP} - \devP^h}_{\Le}$}
    				
    				\addplot[color=asl, mark=square] coordinates {
                            ( 1871 , 0.00016902136346478344 )
                            ( 7223 , 4.190326648221085e-05 )
                            ( 16063 , 7.281008609392417e-06 )
                            ( 21791 , 3.7542369074758815e-06 )
                            ( 28391 , 2.1208602544016252e-06 )
    				};
    				\addlegendentry{$\norm{\widetilde{\sphP} - \sphP^h}_{\Le}$}

                         \addplot[color=violet, mark=pentagon] coordinates {
                            ( 1871 , 0.043680092304945475 )
                            ( 7223 , 0.015812186485471766 )
                            ( 16063 , 0.0026462981100792552 )
                            ( 21791 , 0.00130215837737765 )
                            ( 28391 , 0.0006968336586330491 )
    				};
    				\addlegendentry{$\norm{\widetilde{\hP} - \hP^h}_{\Le}$}
    				
    				\addplot[dashed,color=black, mark=none]
    				coordinates {
    					(6000, 1.5e-4)
    					(20000, 7.394254526319742e-06)
    				};    		
                    \addplot[dashed,color=black, mark=none]
    				coordinates {
    					(7000, 2.5e-5)
    					(20000, 3.0625e-06)
    				};    		
    			\end{loglogaxis}
    			\draw (4.5,1.7) 
    			node[anchor=south west]{$\mathcal{O}(h^{5})$};
                \draw (5.2,1.15) 
    			node[anchor=north east]{$\mathcal{O}(h^{4})$};
    		\end{tikzpicture}
    		\caption{}
    	\end{subfigure}
    	\caption{Convergence rates on the primal (a) and mixed formulations (b) for $\muc = 0$ and $\Lc = 1$.}
    	\label{fig:convex}
\end{figure}
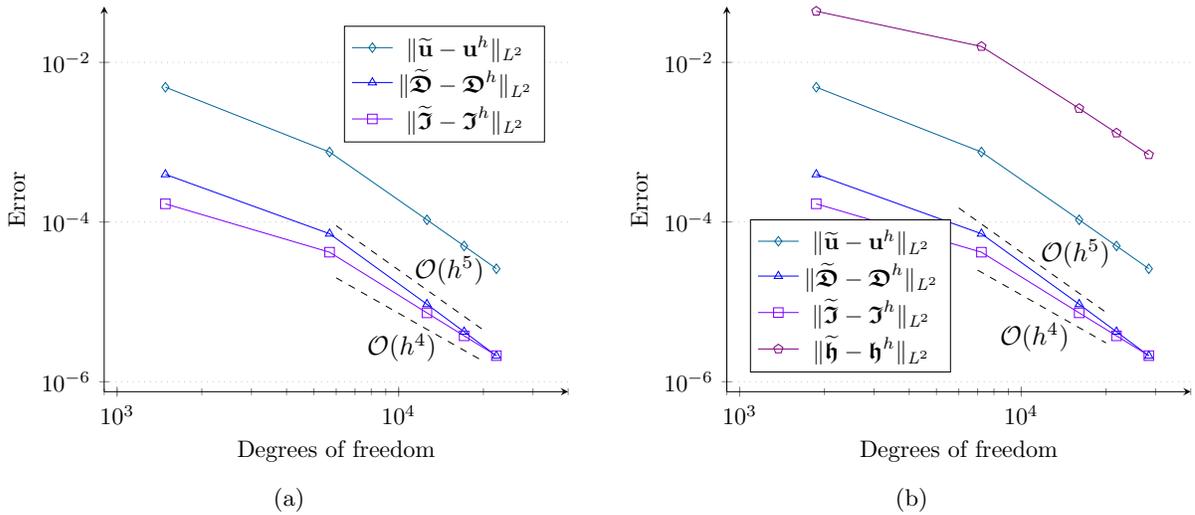

We define the symmetric square domain $\overline{\surf} = [-2\pi,2\pi]^2$ with the homogeneous material coefficients     
\begin{align}
 &\lamma = 115.4 \, , && \muma = 76.9 \, , && \muc = 0 \, , && \Lc = 1 \, , && \lammi = 10 \, \lamma \, , && \mumi = 10 \, \muma \, , 
\end{align}
where the meso-parameters are retrieved via \cref{eq:homogen}. In order to test for convergence we define
\begin{align}
    &\widetilde{\ud} = \dfrac{1}{4} \begin{bmatrix}
        \sin(y) \\
        \sin(x)
    \end{bmatrix} \, , && \widetilde{\devP} = 0 \, , && \widetilde{\sphP} = 0 \, ,
\end{align}
and enforce it on the Dirichlet boundary $\curv_D = \partial \surf$. Consequently, the previous definitions represent the analytical solution if we construct the forces and couple-forces as
\begin{align}
    &\vb{f} = -\Di(\Ce \sym \D \widetilde{\ud}) \, , &&
      \bm{M} = -\Ce  \sym \D \widetilde{\ud} \, . 
\end{align}
The resulting convergence rates are depicted in \cref{fig:convex} for meshes with $32$, $128$, $288$, $392$ and $512$ quartic elements, respectively.
We observe super-optimal convergence rates in all parameters in both the primal and mixed formulations. The improved convergence rates can be explained by the high regularity of the data \cite{KNEES2023126806,Owczarek,knees2023global}, as shown in \cite{SKY2022115298} using $s$-regularity assumptions and the Aubin--Nitsche technique for the full Curl model.   

\begin{figure}
    	\centering
    	\begin{subfigure}{0.48\linewidth}
    		\centering
    		\begin{tikzpicture}[scale = 0.9]
    			\definecolor{asl}{rgb}{0.4980392156862745,0.,1.}
    			\definecolor{asb}{rgb}{0.,0.4,0.6}
    			\begin{loglogaxis}[
    				/pgf/number format/1000 sep={},
    				axis lines = left,
    				xlabel={Degrees of freedom},
    				ylabel={Error},
    				xmin=900, xmax=40000,
    				ymin=0.75e-6, ymax=0.05,
    				xtick={1e2,1e3,1e4,1e5},
    				ytick={1e-6,1e-4,1e-2},
    				legend style={at={(0.95,0.95)},anchor= north east},
    				ymajorgrids=true,
    				grid style=dotted,
    				]
    				\addplot[color=asb, mark=diamond] coordinates {
    					( 1485 , 0.029361319856003567 )
                            ( 5685 , 535.5489954772761 )
                            ( 12605 , 0.00010850329765298439 )
                            ( 17085 , 5.7325224705796624e-05 )
                            ( 22245 , 2.63436162984537e-05 )
    				};
    				\addlegendentry{$\norm{\widetilde{\ud} - \ud^h}_{\Le}$}

                    \addplot[color=blue, mark=triangle] coordinates {
    					( 1485 , 2.5057081169181877e-06 )
                            ( 5685 , 0.003662360106465721 )
                            ( 12605 , 4.9859547579161123e-11 )
                            ( 17085 , 3.77695001476862e-11 )
                            ( 22245 , 1.0589465050762703e-12 )
    				};
    				\addlegendentry{$\norm{\widetilde{\devP} - \devP^h}_{\Le}$}
    				
    				\addplot[color=asl, mark=square] coordinates {
                            ( 1485 , 0.003682391566864882 )
                            ( 5685 , 102.53486004843955 )
                            ( 12605 , 8.029519995051916e-06 )
                            ( 17085 , 9.322293398982352e-06 )
                            ( 22245 , 2.1644682522473732e-06 )
    				};
    				\addlegendentry{$\norm{\widetilde{\sphP} - \sphP^h}_{\Le}$}
    				
    				\addplot[dashed,color=black, mark=none]
    				coordinates {
    					(6000, 9e-5)
    					(20000, 4.436552715791846e-06)
    				};    		
                    \addplot[dashed,color=black, mark=none]
    				coordinates {
    					(6000, 2e-5)
    					(20000, 1.8e-6)
    				};    		
    			\end{loglogaxis}
    			\draw (4.45,1.45) 
    			node[anchor=south west]{$\mathcal{O}(h^{5})$};
                \draw (5.05,1) 
    			node[anchor=north east]{$\mathcal{O}(h^{4})$};
    		\end{tikzpicture}
    		\caption{}
    	\end{subfigure}
    	\begin{subfigure}{0.48\linewidth}
    		\centering
    		\begin{tikzpicture}[scale = 0.9]
    			\definecolor{asl}{rgb}{0.4980392156862745,0.,1.}
    			\definecolor{asb}{rgb}{0.,0.4,0.6}
    			\begin{loglogaxis}[
    				/pgf/number format/1000 sep={},
    				axis lines = left,
    				xlabel={Degrees of freedom},
    				ylabel={Error},
    				xmin=900, xmax=40000,
    				ymin=0.75e-6, ymax=0.05,
    				xtick={1e2,1e3,1e4,1e5},
    				ytick={1e-6,1e-4,1e-2},
    				legend style={at={(0.05,0.05)},anchor= south west},
    				ymajorgrids=true,
    				grid style=dotted,
    				]
    				\addplot[color=asb, mark=diamond] coordinates {
                            ( 1871 , 0.004893783417188508 )
                            ( 7223 , 0.0007588201817025512 )
                            ( 16063 , 0.00010779085426705161 )
                            ( 21791 , 5.074224094454129e-05 )
                            ( 28391 , 2.63297347407172e-05 )
    				};
    				\addlegendentry{$\norm{\widetilde{\ud} - \ud^h}_{\Le}$}

                    \addplot[color=blue, mark=triangle] coordinates {
                            ( 1871 , 1.2962736867666764e-05 )
                            ( 7223 , 6.900668858977677e-08 )
                            ( 16063 , 3.874833591805034e-09 )
                            ( 21791 , 1.3017615387432133e-09 )
                            ( 28391 , 5.0716794206018e-10 )
    				};
    				\addlegendentry{$\norm{\widetilde{\devP} - \devP^h}_{\Le}$}
    				
    				\addplot[color=asl, mark=square] coordinates {
                            ( 1871 , 0.00017623420215662 )
                            ( 7223 , 4.307792648759992e-05 )
                            ( 16063 , 7.438779006003042e-06 )
                            ( 21791 , 3.826628694367072e-06 )
                            ( 28391 , 2.1577115203148773e-06 )
    				};
    				\addlegendentry{$\norm{\widetilde{\sphP} - \sphP^h}_{\Le}$}

                         \addplot[color=violet, mark=pentagon] coordinates {
                            ( 1871 , 1.9071071073135848 )
                            ( 7223 , 0.16497488836675783 )
                            ( 16063 , 0.04042994337339441 )
                            ( 21791 , 0.011723957305014052 )
                            ( 28391 , 0.010795645481568486 )
                            ( 44207 , 0.0032459331170882825 )
    				};
    				\addlegendentry{$\norm{\widetilde{\hP} - \hP^h}_{\Le}$}
    				
    				\addplot[dashed,color=black, mark=none]
    				coordinates {
    					(6000, 1.5e-4)
    					(20000, 7.394254526319742e-06)
    				};    		
                    \addplot[dashed,color=black, mark=none]
    				coordinates {
    					(7000, 2.5e-5)
    					(20000, 3.0625e-06)
    				};    		
    			\end{loglogaxis}
    			\draw (4.5,1.7) 
    			node[anchor=south west]{$\mathcal{O}(h^{5})$};
                \draw (5.2,1.15) 
    			node[anchor=north east]{$\mathcal{O}(h^{4})$};
    		\end{tikzpicture}
    		\caption{}
    	\end{subfigure}
    	\caption{Convergence rates on the primal (a) and mixed formulations (b) for $\muc = 1$ and $\Lc = 10^{10}$.}
    	\label{fig:convexlc}
\end{figure}
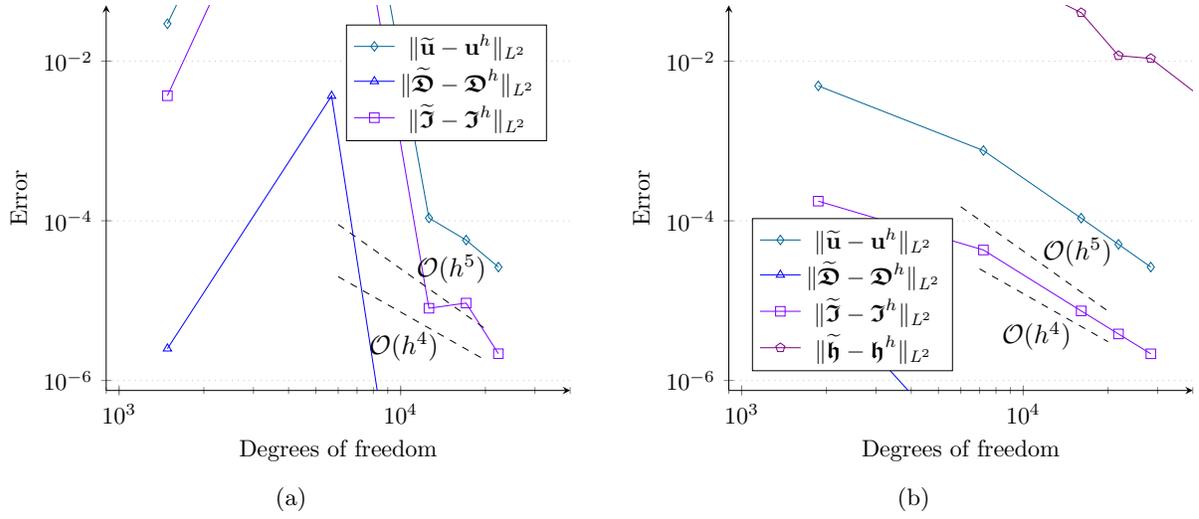

Considering the same example but with the parameters
\begin{align}
    &\muc = 1 \, , && \Lc = 10^{10} \, ,
\end{align}
representing the case $\Lc \to +\infty$ of a zoom into the micro-structure, we observe vastly different results, as depicted in \cref{fig:convexlc}.  
Namely, the primal formulation becomes unstable and the quality of the solution oscillates. In contrast, the mixed formulation remains optimal in the displacement $\ud$, spherical microdistortion $\sphP$ and the hyper-stress $\hP$, while the deviatoric microdistortion $\devP$ converges even faster. We note that in comparison, for $\muc = 0$ and $\Lc = 10^{10}$ the finite element computation is more stable.

\begin{figure}
    	\centering
    	\begin{subfigure}{0.48\linewidth}
    		\centering
    		\definecolor{npurple}{rgb}{0.4980392156862745,0.,1.}
\begin{tikzpicture}[scale = 0.9]
			\begin{semilogxaxis}[
				/pgf/number format/1000 sep={},
				axis lines = left,
				xlabel={$\Lc$},
				ylabel={Energy},
				xmin=0.5e-3, xmax=5e+3,
				ymin=100, ymax=120,
				x dir=reverse,
				xtick={1e-4, 1e-2, 1, 1e+2, 1e+4},
				ytick={105,110,115},
				legend pos= north east,
				ymajorgrids=true,
				grid style=dotted,
				]
				\addplot[
				color=npurple,
				mark=triangle,
				]
				coordinates {
					( 1000.0 , 113.9685419417981 )
                        ( 100.0 , 113.92096109768129 )
                        ( 31.622776601683793 , 113.5292729631851 )
                        ( 17.78279410038923 , 112.78120968475953 )
                        ( 14.12537544622754 , 112.28067426219044 )
                        ( 10.0 , 111.30553284432173 )
                        ( 5.623413251903491 , 109.39198116194783 )
                        ( 3.1622776601683795 , 107.71090378171272 )
                        ( 1.2589254117941673 , 106.28098009248357 )
                        ( 1.0 , 106.10663851828185 )
                        ( 0.31622776601683794 , 105.59716594081993 )
                        ( 0.1 , 105.32881428238609 )
                        ( 0.01 , 105.17900677442827 )
                        ( 0.001 , 105.16334463953424 )
				};
				\addlegendentry{$I(\ud,\devP,\sphP)$}
				\addplot[dashed,color=black, mark=none]
				coordinates {
					(1e+4, 600)
					(1e-4, 600)
				};
				\addplot[dashdotted,color=black, mark=none]
				coordinates {
					(1e+4, 514.0757380591552)
					(1e-4, 514.0757380591552)
				};
				\addplot[dashed,color=black, mark=none]
				coordinates {
					(1e+4, 105.0160649870426)
					(1e-4, 105.0160649870426)
				};
			\end{semilogxaxis}
			\draw (1.,1.4) node[anchor=north west]{$\Cmac$};
		\end{tikzpicture}
    		\caption{}
    	\end{subfigure}
    	\begin{subfigure}{0.48\linewidth}
    		\centering
    		\includegraphics[width=0.8\linewidth]{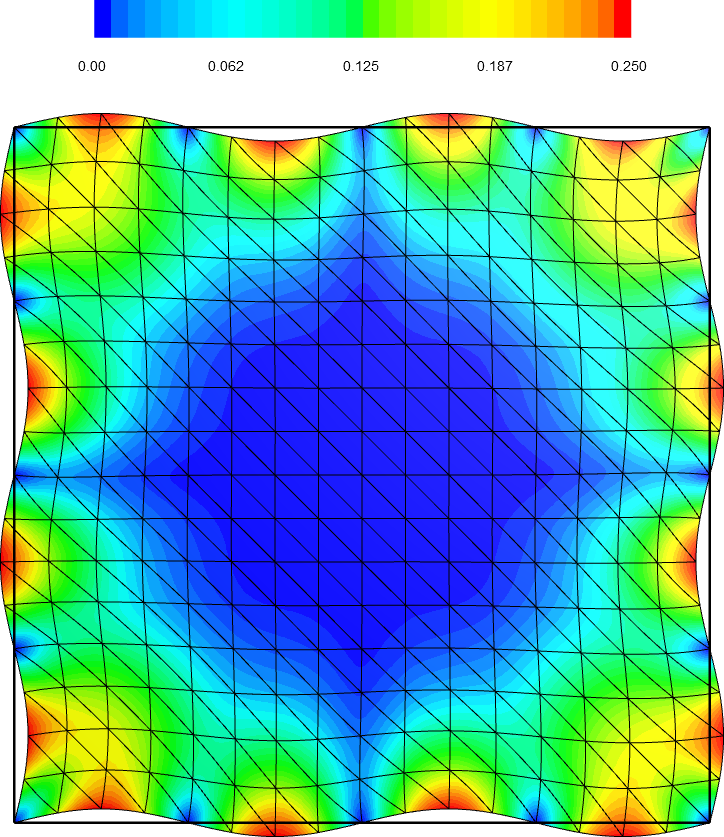}
    		\caption{}
    	\end{subfigure}
    	\caption{Convergence in the elastic energy of the relaxed micromorphic model for $\Lc \to 0$ towards the equivalent energy of linear elasticity with $\Cmac$ (a), and the corresponding displacement field $\ud$ for $\Lc = 10^{-3}$ (b).}
    \label{fig:energylc}
\end{figure}
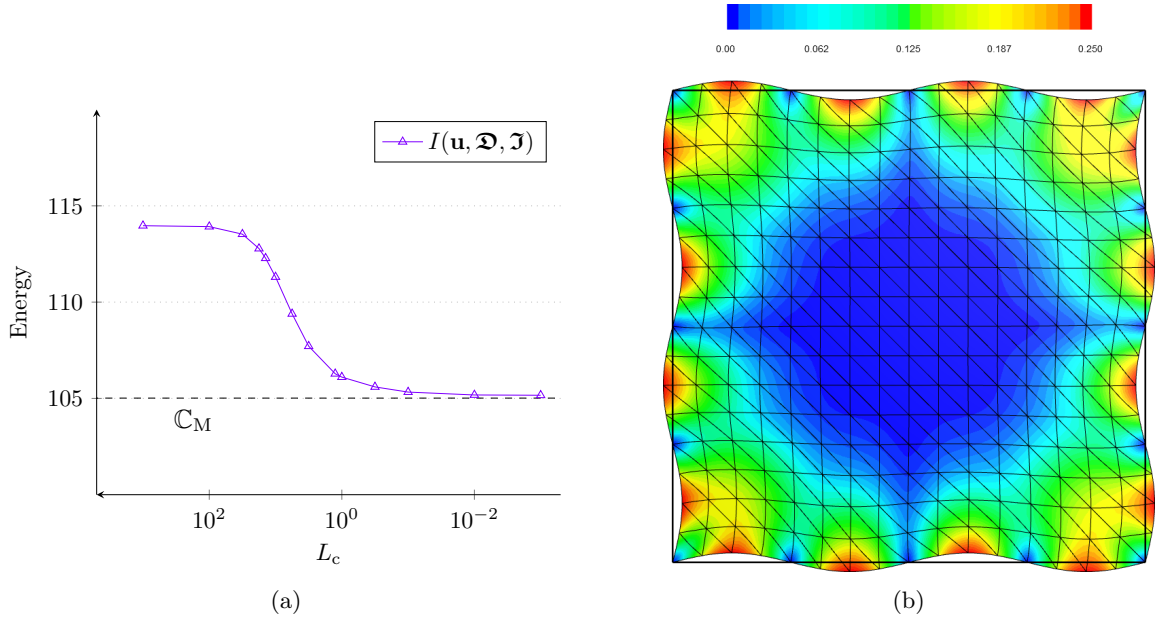

Finally, setting the forces and couple-forces to zero $\vb{f} = \bm{M} = 0$, with $\muc = 0$ and letting $\Lc \to 0$, we expect to retrieve an equivalent energy $I$ to that of a homogeneous  linear elastic model with $\Cmac$, as per \cref{eq:lctozero}.  
We employ the fine mesh of $512$ quartic elements and vary the characteristic length-scale parameter $\Lc \in [10^3,10^{-3}]$. The lower limit is given by a computation of linear elasticity with $\Cmac$. The change in elastic energy of the relaxed micromorphic model is depicted in \cref{fig:energylc}. 
We can clearly observe convergence towards the energy given by the macro moduli, such that the energies practically agree for $\Lc = 10^{-3}$. Further, for $\Lc = 10^3$ the energy stagnates as well, implying a finite limit. However, as shown in \cref{sec:lclimits} we cannot characterise this limit as micro-level elasticity problem for the $\sym\Curl$-model.

\section{Conclusions and outlook}
In this work we chose to split the microdistortion field of the relaxed micromorphic model into its deviatoric and spherical parts, explicitly. Doing so allowed us to consistently reduce the model to two dimensions without losing its capacity to capture jumping dilatation fields, which may occur in micro-structured materials. We introduced the reduced model and explored its behaviour for the limits of the characteristic length-scale parameter $\Lc \in [0,+\infty)$, leading to primal and mixed variational formulations, as well as novel micro-scale differential equations.
The isochoric split of the microdistortion field and subsequent dimensional reduction required the introduction of new conforming finite element spaces for $\HCd{,\surf}$. We proposed two such elements and proved their linear independence, or unisolvence, respectively. The viability of our approach as well as the correctness of the first finite element space were demonstrated in two numerical benchmarks. In the first example, jumping material coefficients inherently led to a jump in the dilatation, while in the second example the stability of the mixed formulation for $\Lc \to +\infty$ was shown. 

The implementation of the macro element approach is not trivial and is thus reserved for future works. In cases where a jump in the dilatation field and the deviatoric field occur simultaneously, this element becomes indispensable, as it can correctly account for both, the minimal $\HCd{,\surf}$-regularity, and consistent boundary conditions.  

\section*{Acknowledgements}
Patrizio Neff acknowledges support in the framework of the DFG-Priority Programme
2256 “Variational Methods for Predicting Complex Phenomena in Engineering Structures and Materials”, Neff 902/10-1, Project-No. 440935806. Michael Neunteufel acknowledges support by the US National Science Foundation.

\bibliographystyle{spmpsci}   

\footnotesize
\bibliography{Ref}   

\normalsize
\appendix
\section{Two-dimensional curl operators} \label{ap:b}
In the three-dimensional Cartesian system the curl operator reads
\begin{align}
    \curl \vb{v} = \nabla \times \vb{v} = \begin{bmatrix}
        v_{3,y} - v_{2,z} \\
        v_{1,z} - v_{3,x} \\
        v_{2,x} - v_{1,y} 
    \end{bmatrix} \, .
\end{align}
Its row-wise application to second order tensors is therefore 
\begin{align}
    \Curl \Pm = \begin{bmatrix}
        (\curl \vb{p}_1)^T \\
        (\curl \vb{p}_2)^T \\
        (\curl \vb{p}_3)^T 
    \end{bmatrix} = - \Pm \times \nabla \, ,
\end{align}
where $\vb{p}_i$ represent rows of $\Pm$.
Now, if the vector field $\vb{v}$ is planar and depends only on the variables $x$ and $y$, then its curl reduces to
\begin{align}
    \curl \vb{v} = \nabla \times \vb{v} = \begin{bmatrix}
        0 \\
        0 \\
        v_{2,x} - v_{1,y} 
    \end{bmatrix}  \, .
\end{align}
Assuming a two-dimensional view of the problem, fields in $\vb{e}_3$ are interpreted as scalar intensities and we find
\begin{align}
    &\curl \vb{v} = \di (\bm{R} \vb{v}) = v_{2,x} - v_{1,y} \,, && \bm{R} = \begin{bmatrix}
        0 & 1 \\
        -1 & 0
    \end{bmatrix} \, , && \vb{v} = \begin{bmatrix}
        v_1 \\
        v_2
    \end{bmatrix} \, .
\end{align}
Alternatively, if $v_1 = v_2 = 0$ and $v_3 = v_3(x,y)$ we find
\begin{align}
    \curl \vb{v} = \nabla \times \vb{v} = \begin{bmatrix}
        v_{3,y} \\
        -v_{3,x} \\
        0
    \end{bmatrix} \, .
\end{align}
Since components in $\vb{e}_1$ and $\vb{e}_2$ are interpreted as vectorial fields of the plane, we can recast the operation as
\begin{align}
    \bm{R} \nabla \lambda = \begin{bmatrix}
        \lambda_{y} \\
        -\lambda_{x} 
    \end{bmatrix} \, ,
\end{align}
where $v_3 = \lambda$ is interpreted as a scalar field of the problem. Thus, in two dimensions we find two operators representing the curl of a field, one for vectors $\curl(\cdot) = \di(\bm{R}\cdot)$ and one for scalars $\bm{R}\nabla(\cdot)$. 
Consequently, the application of the vectorial curl operator to two-dimensional second order tensors reads
\begin{align}
    &\Curl \Pm = \begin{bmatrix}
        \curl \vb{p}_1 \\
        \curl \vb{p}_2 
    \end{bmatrix}
     = \begin{bmatrix}
        \di (\bm{R} \vb{p}_1) \\
        \di (\bm{R} \vb{p}_2) 
    \end{bmatrix}= \Di (\Pm\bm{R}^T)  = -(\Pm \cdot \bm{R}) \cdot \nabla \, , 
\end{align}
whereas the scalar curl operator is applied to vectors 
\begin{align}
    \begin{bmatrix}
         (\bm{R} \nabla v_1)^T \\
         (\bm{R} \nabla v_2)^T \\
    \end{bmatrix} = (\D \vb{v}) \bm{R}^T = -(\vb{v} \otimes \nabla) \bm{R} \, .
\end{align}

\section{The $\sym\Curl$ and $\Curl$ operators and corresponding traces}\label{ap:a}
Spherical fields are in the kernel of the $\sym\Curl$-operator
\begin{align}
    \sym\Curl\Pm &= \sym\Curl(\dev\Pm + \sph \Pm) 
    = \sym\Curl \dev\Pm + \dfrac{1}{3}\sym\Curl[(\tr\Pm) \one]
     \\
    &= \sym\Curl \dev\Pm - \dfrac{1}{3}\sym\Anti(\nabla\tr\Pm) = \sym\Curl \dev\Pm  \, , \notag
\end{align}
where $\Anti(\cdot)$ yields the skew-symmetric matrix representation of the cross product.
Similarly, for the trace operator we observe 
\begin{align}
    \tr_{\HsC{}}\Pm &= \sym[\Pm (\Anti \vb{n})^T] \at_{\Xi} 
    = \sym[(\dev\Pm + \sph \Pm) (\Anti \vb{n})^T] \at_{\Xi} 
     \\
    &= \sym[(\dev\Pm) (\Anti \vb{n})^T] + \dfrac{1}{3}\sym[ (\tr\Pm) \one (\Anti \vb{n})^T] \at_{\Xi}
    = \sym[(\dev\Pm) (\Anti \vb{n})^T] \at_{\Xi} \,, \notag
\end{align}
where $\Xi \subset \R^2$ represents some arbitrary interface or boundary and $\vb{n}$ is its normal vector. It is directly observable that without the $\sym$-operator the spherical part is not eliminated in the differentiation  
\begin{align}
    \Curl\Pm &= \Curl(\dev\Pm + \sph \Pm) 
    = \Curl \dev\Pm - \dfrac{1}{3}\Anti(\nabla\tr\Pm) \, ,
\end{align}
or in the trace
\begin{align}
    \tr_{\HC{}}\Pm 
    &= (\dev\Pm + \sph \Pm) (\Anti \vb{n})^T \at_{\Xi} 
    = (\dev\Pm) (\Anti \vb{n})^T + \dfrac{1}{3}(\tr\Pm)  (\Anti \vb{n})^T \at_{\Xi} \, .
\end{align}
Consequently, the decomposition
\begin{align}
    \HsC{,\body} = \HsCd{,\body} \oplus \HsCv{,\body}  = \HsCd{,\body} \oplus [\Le(\body) \cdot \one] \, ,
\end{align}
does not correspond to
\begin{align}
    \HC{,\body} = \HCd{,\body} \oplus \HCv{,\body} \, ,
\end{align}
since 
\begin{align}
    \HCv{,\body} = \{ \Pm \in [\Le(\body) \cdot \one] \; | \; \nabla \tr \Pm  \in [\Le(\body)]^3 \} = [\Hone(\body) \cdot \one] \neq  [\Le(\body) \cdot \one] \, .
\end{align}
Analogously, in two dimensions we have
\begin{align}
    \Curl \Pm &= \Curl\dev \Pm + \Curl\sph\Pm 
    = \Di[\dev \Pm \bm{R}^T] + \dfrac{1}{2} \Di[(\tr \Pm) \one \bm{R}^T] = \Di[\dev \Pm \bm{R}^T] + \dfrac{1}{2} \bm{R}^T \nabla(\tr \Pm) \, ,
\end{align}
for the differential operator, and
\begin{align}
    \tr_{\HC{}}\Pm &= (\dev \Pm + \sph\Pm)\vb{t} \at_\Xi 
    = (\dev \Pm) \vb{t} + \dfrac{1}{2} (\tr \Pm) \one \vb{t} \at_\Xi  
    = (\dev \Pm) \vb{t} + \dfrac{1}{2} (\tr \Pm)  \vb{t} \at_\Xi  \, ,
\end{align}
for the trace operator, where $\vb{t}$ is the tangent vector on $\Xi \subset \R$. 

\end{document}